\documentclass{article}
\usepackage{amsmath,amssymb,amscd,amsthm,amsfonts,amstext,amsxtra}
\usepackage{eufrak,fontenc,fontsmpl,latexsym,rawfonts}
\theoremstyle{plain}
\newtheorem{theorem}{Theorem}
\newtheorem{lemma}[theorem]{Lemma}

\newtheorem{definition}[theorem]{Definition}

\newtheorem*{theorem*}{Theorem}
\newtheorem*{lemma*}{Lemma}
\newtheorem*{proposition*}{Proposition}
\newtheorem*{corollary*}{Corollary}
\newtheorem*{definition*}{Definition}
\newtheorem*{conjecture*}{Conjecture}
\theoremstyle{definition}
\newtheorem*{remark*}{Remark}
\newtheorem*{acknowledgements}{Acknowledgements}
\newtheorem*{keywords}{Keywords}
\newtheorem{remark}[theorem]{Remark}
\begin{document}
\title{Direct image for multiplicative and relative
$K$-theories from transgression of the families index theorem, part 1.}
\author{Alain Berthomieu\\\and Institut de Math\'ematiques de Toulouse, UMR CNRS n$^\circ$ 5219\and et Centre Universitaire de Formation et de Recherche\and Jean-Fran\c{c}ois Champollion, Campus d'Albi,\and Place de Verdun,
81012 Albi Cedex, France.}
\maketitle
\begin{abstract}
This paper is the first part of the ``longstanding forthcoming preprint'' referred to in \cite{MoiOberwolfach}. It contains the constructions
of the real counterpart of the relative $K$-theory considered in \cite{BerthomieuPNCDIRKT} in the context of complex analytic geometry,
and of an extension of Karoubi's multiplicative $K$-theory \cite{KaroubiCCFFHA} suggested by U. Bunke \cite{BunkeSchick}
(which I call ``free multiplicative $K$-theory'' in the sequel).

Nadel's type classes on relative $K$-theory as in \cite{Nadel} and \cite{BerthomieuPNCDIRKT}
are constructed, while it is proved that on free multiplicative $K$-theory,
there is a notion of Chern-Weil character form, and of a Borel-type
characteristic class (which is a differential form modulo exact forms)
which recovers the classes $c_k$ studied by Bismut and Lott in \cite{BismutLott}.

Finally, a direct image for relative $K$-theory under proper submersion
of compact orientable real manifolds is constructed.
\begin{keywords} multiplicative $K$-theory, relative $K$-theory, Chern-Simons transgression, superconnections, flat vector bundles,
proper submersions, direct image.
\end{keywords}
\noindent{\textbf{AMS-classification:}} Primary: 14F05, 19E20, 57R20, secondary: 14F40, 19D55, 53C05, 55R50.
\end{abstract}
\section{Introduction:}

Consider four complex vector bundles $E$, $F$, $G$ and $H$ on some
manifold $M$. If it happens that $[E]-[F]=[G]-[H]\in K^0_{\text{top}}(M)$
then there exists some vector bundle $K$ and some smooth isomorphism
$f\colon E\oplus H\oplus K\overset\sim\longrightarrow F\oplus G\oplus K$.
Such data $(K,f)$ will be called a ``link'' between $E-F$ and $G-H$.
In fact a link can be twisted by some element of $K^1_{\text{top}}(M)$.
But there exists an equivalence relation between links which measures
such twists (see \S\ref{S:linksetrelations}).

Let $\pi\colon M\to B$ be a proper smooth submersion between compact manifolds,
explicit constructions of direct image morphism
$\pi_*^{\rm{Eu}}\colon K^0_{\text{top}}(M)\to K^0_{\text{top}}(B)$ by analytic families index theory (using the fibral Euler ($\cong$signature) operator) as in \cite{AtiyahSinger},
\cite{MischchenkoFomenko} and \S\ref{directimtopologic} produce, for some vector bundle
$\xi$ on $M$, some couples of vector bundles on $B$ whose differences represent
$\pi_*^{\rm{Eu}}[\xi]\in K^0_{\text{top}}(B)$, and it is proved in paragraph \ref{clink} below
that all these couples can be canonically linked in the above sense.

The first goal of this paper is to construct a real analogue $K^0_{\text{rel}}$
of the relative
$K$-theory of \cite{BerthomieuPNCDIRKT} whose objects are quintuples $(E,\nabla_{\!E}, F,\nabla_{\!F},f)$
where $E$ and $F$ are complex vector bundles, $\nabla_{\!E}$ and $\nabla_{\!F}$
are flat connections on them, and $f\colon E\to F$ is a
smooth isomorphism. The above linking of representatives of direct image in
topological $K$-theory is then needed to construct the direct image morphism
$\pi_*\colon K^0_{\text{rel}}(M)\to K^0_{\text{rel}}(B)$ in \S\ref{dirimrelsection} below.

In such a situation, characteristic classes are constructed
using Chern-Simons forms: let ${\rm{ch}}(\nabla_{\!E})$ and ${\rm{ch}}(\nabla_{\!F})$ be the
explicit differential form representatives of the Chern
characters of $E$ and $F$ obtained throw Chern-Weil theory
from the connections $\nabla_{\!E}$ and $\nabla_{\!F}$, then
Chern-Simons theory provides an odd differential form
$\widetilde{\rm{ch}}(\nabla_{\!E},f^*\nabla_{\!F})$ defined modulo coboundaries which verifies
\[d\widetilde{\rm{ch}}(\nabla_{\!E},f^*\nabla_{\!F})={\rm{ch}}
(\nabla_{\!F})-{\rm{ch}}(\nabla_{\!E})\]
(This of course also works for two nonnecessary flat connexions on any same vector bundle). Here the curvatures of $\nabla_{\!E}$ and $\nabla_{\!F}$ vanish so
that $\widetilde{\rm{ch}}(\nabla_{\!E},f^*\nabla_{\!F})$ is a closed form, and one gets a group morphism ${\mathcal N}_{\rm{ch}}\colon K_{\text{rel}}^0(M)
\longrightarrow H^{\text{odd}}(M,{\mathbb C})$.

Chern-Simons forms can also be used to construct a ``free multiplicative'' $K$-theory
$\widehat K_{\rm{ch}}(M)$ whose objects are triples $(E,\nabla,\alpha)$ where $E$ is a complex vector
bundle on $M$ with connection $\nabla$ and $\alpha$ is an odd-degree differential form on $M$
defined modulo exact forms, with following relations: direct sum, and
change of connection in the following way
$(E,\nabla_{\!0},\alpha)=\big(E,\nabla_{\!1},\alpha+\widetilde{\rm{ch}}(\nabla_{\!0},\nabla_{\!1})\big)$.
On such an object, one has a Chern character differential form $\widehat{\rm{ch}}(E,\nabla,\alpha)=
{\rm{ch}}(\nabla)-d\alpha$
and a Borel-type class ${\mathfrak{B}}(E,\nabla,\alpha)=
\widetilde{\rm{ch}}(\nabla^*,\nabla)-2i\, {\mathfrak{Im}}\, \alpha$
(this is a purely imaginary odd differential form defined modulo exact forms)
where $\nabla^*$ is the adjoint transpose of the connection
$\nabla$ with respect to any hermitian metric on $E$. This free multiplicative
$K$-theory with differential form Chern character generalises Karoubi's multiplicative $K$-theory \cite{KaroubiHCKT} \cite{KaroubiTGCCS}
\cite{KaroubiCCFFHA} and
provides classes for vector bundles with connections.

The paper is organized as follows:
in section \ref{various}, the $K$-theories considered are
recalled or constructed, and the considered caracteristic
classes on them are introduced, in section \ref{directim},
the direct image for topological and relative $K$-theories under
smoth submersions are
recalled or constructed.

I plan to construct a direct image in the same context for
free multiplicative $K$-theory in a second part which would be consistent with the usual sheaf theoretic direct image of flat vector bundles. As a by-product of this construction, a
Riemann-Roch-Grothendieck type theorem for Nadel's classes on
relative $K$-theory should be obtained.

Finally, the question of double fibrations should be studied in a third part.
\begin{acknowledgements}
I am much indebted to Ullrich Bunke for suggesting to relax the condition
\eqref{DankeBunke} in the definition of transgressive $K$-theory \cite{BunkeSchick}, and to the organizers of the workshop ``Analysis and topology in interaction'' held in Oberwolfach in March 2006 for their invitation and the nice and stimulating discussions which took place there.
\end{acknowledgements}
\section{Various $K$-theories:}\label{various}
\subsection{Topological, flat and relative $K$-theories:}
\subsubsection{Topological $K$-theory:}\label{cathop}
Let $M$ be any ${\mathcal C}^\infty$ real manifold.
As usual, the topological $K$-group $K^0_{\text{top}}(M)$
is the free abelian group generated by isomorphism classes
of smooth complex vector bundles on $M$ modulo direct sum.

The topological $K^1$-group $K^1_{\text{top}}(M)$ can either be
seen as the quotient group $K^0_{\text{top}}(M\times S^1)/p_1^*
K^0_{\text{top}}(M)$ where $p_1\colon M\times S^1\to M$ is
the obvious projection, or as the kernel of the
restriction map $\iota^*\colon K^0_{\text{top}}(M\times S^1)\to
K^0_{\text{top}}(M\times\{pt\})$ where $pt$ is some point in
$S^1$ and $\iota\colon pt\to S^1$ the inclusion map.

Another description we will use, is that $K^1_{\text{top}}(M)$
is generated by global smooth automorphisms $g$ of trivial
complex vector bundles ${\mathbb C}^N$ on $M$ (of all possible
ranks $N$). The corresponding vector bundle on $M\times S^1$
is obtained by gluing using $g$ the restrictions to
$M\times\{1\}$ and $M\times \{0\}$ of the trivial vector bundle
${\mathbb C}^N$ on $M\times [0,1]$. An element of $K^1_{\rm{top}}(M)$
can be described as a global smooth vector bundle automorphism $g_E$ of any vector bundle $E$ on $M$, but the corresponding vector bundle on
$M\times S^1$ should be constructed by first adding to $E$ another vector bundle $F$ such that $E\oplus F$ is trivial, and glue the two copies of $E\oplus F$
over $M\times\{0\}$ and $M\times\{1\}$ throw $g_E\oplus{\rm{Id}}_F$.

If $g_E={\rm{Id}}_E$, the correponding class in $K^1_{\rm{top}}(M)$ is trivial,
and this remains true if $g_E$ is a constant multiple of the identity. If
$g_E$ is a selfadjoint positive (with respect to any hermitian metric)
automorphism of $E$, the associated class in $K^1_{\rm{top}}(M)$ is also
trivial, because the space of selfadjoint positive automorphisms (with respect
to the same hermitian metric) is convex and contains the identity. If
$g_E$ is selfadjoint, $E$ decomposes as an orthogonal direct sum $E=E^>\oplus E^<$ of
subvector bundles respected by $g_E$ on which $g_E$ is selfadjoint positive and
selfadjoint negative respectively. The element of $K^1_{\rm{top}}(M)$ associated to each resriction of $g_E$ vanish, so that the one associated to
$g_E$ also does. Finally, an argument of the same sort works to prove the triviality of the element of $K^1_{\rm{top}}(M)$ provided by a skew-adjoint
automorphism of $E$ (with respect to any hermitian metric on $E$).
\subsubsection{Links and complexes of vector bundles:}\label{S:linksetrelations}

For four smooth vector bundles $E$, $F$, $G$, and $H$ on $M$
such that
\[[E]-[F]=[G]-[H]\in K^0_{\text{top}}(M)\]
there
exists some vector bundle $K$ on $M$ and some ${\mathcal C}^\infty$ isomorphism
\begin{equation}\label{linque}\ell\colon E\oplus H\oplus K\overset\sim\longrightarrow F\oplus G\oplus K\end{equation}
These $(K,\ell)$ will be called a ``link between
$E-F$ and $G-H$''. Two such links $(K_1,\ell_1)$ and $(K_2,\ell_2)$ are
equivalent if there exists some vector bundle $L$ on $M$ such that the two isomorphisms
\begin{equation}\label{linkdelink}E\oplus H\oplus K_1\oplus K_2\oplus L\overset
{\ell_1\oplus{\text{Id}}_{K_2}\oplus{\text{Id}}_{L}}{\underset
{\ell_2\oplus
{\text{Id}}_{K_1}\oplus{\text{Id}}_L}{-\!\!\!-\!\!\!-\!\!\!-\!\!\!-\!\!\!-\!\!\!-\!\!\!\-\!\!\!-\!\!\!-\!\!\!-\!\!\!\longrightarrow}}
F\oplus G\oplus K_1\oplus K_2\oplus L\end{equation}
are isotopic, which here means that they are homotopic
throw isomorphisms. Any link is equivalent to some other one with a trivial vector bundle as $K$. The equivalence class of a link $(K,\ell)$ will be denoted by $[\ell]$.

The expression ``link between\ldots'' is not
commutative: if $(K,\ell)$ is a link between $E-F$ and $G-H$, then
$(K,\ell^{-1})$ will be a link between $G-H$ and $E-F$. Its equivalence class will be denoted by $[\ell^{-1}]$ or $[\ell]^{-1}$.

Links
can be pulled back, and added (for direct sum of data). Moreover, two links $(L,\ell)$ between $E-F$ and $G-H$, and $(M,\ell')$
between $G-H$ and $J-K$ can be composed as
$(L\oplus M\oplus G\oplus H,\ell\oplus\ell')$
between $E-F$ and $J-K$;
this composition is associative. The equivalence class of the composed link will be denoted by $[\ell'\circ\ell]$ or $[\ell']\circ[\ell]$.

If $[\ell]$ is an equivalence class of link between $E-F$ and $G-H$ the obviously associated link between $F-E$ and $H-G$
will be denoted $[-\ell]$. 

An exact sequence of complex vector bundles on $M$ as
\begin{equation}\label{precedente}
0\to E\longrightarrow G\longrightarrow H\longrightarrow F\to 0
\end{equation}
defines unambiguously an equivalence class of
link from $E-F$ to
$G-H$, which is not modified by isotopy of the exact sequence
(``isotopy'' here means ${\mathcal C}^\infty$ homotopy of the
morphisms such that the sequence remains exact at each stage). 

In the same way, four terms or six-terms exact sequence as
\begin{equation}\label{complexestousrondsavant}
\begin{CD}E@>>>G\\@AAA@VVV\\F@<<<H\end{CD}\qquad{\text{ or }}\qquad
\begin{CD}A@>>>B@>>>C\\@AAA@.@VVV\\F@<<<E@<<<D\end{CD}
\end{equation}
produces an unambiguous isotopy-invariant equivalence class of link between $E-F$ and $G-H$ or between $B-E$
and $(A\oplus C)-(D\oplus G)$. Note however that the ranks of the morphisms
may vary along the manifold in contrast with the situation of \eqref{precedente} above.

This can be generalized for a complex of vector bundles on $M$ as
\begin{equation}\label{E:linkcohom}0\overset{v_0}\longrightarrow E^1\overset{v_1}\longrightarrow E^2\overset{v_2}\longrightarrow\ldots\overset{v_{i-1}}\longrightarrow E^i
\overset{v_i}\longrightarrow\ldots\overset{v_{k-1}}\longrightarrow E^k\overset{v_k}\longrightarrow0\end{equation}
where of course $v_i\circ v_{i-1}$ vanishes for any $1\leq i\leq k$.
Suppose that for any $1\leq i\leq k$, the cohomology of $v$ in degree $i$, namely
${\rm{Ker}}v_i\big/{\rm{Im}}v_{i-1}$ is a vector bundle on
$M$, and call it $H^i$. Put $E^+=\big(\mathop{\oplus}\limits_{i\ {\text{even}}}E^i\big)$, $E^-=\big(\mathop{\oplus}\limits_{i\ {\text{odd}}}E^i\big)$ and define $H^+$ and $H^-$ similarly.
Then there exists a canonical equivalence class of link
between $E^+-E^-$ and $H^+-H^-$ which is constructed
in the following way: choose
bundle maps $w_i\colon E^i\to E^{i-1}$ such that $v_{i-1}\circ w_i$ induces
the identity on the image of $v_{i-1}$ in $E^i$,
choose bundle maps $p_i\colon E^i\to H^i$
inducing the canonical projections of ${\rm{Ker}}v_i$ onto $H^i$ and choose bundle maps
$h_i\colon H^i\to E^i$, with images included in ${\rm{Ker}}v_i\subset E^i$, such that $p_i\circ h_i$ is the identuty of $H^i$ for all $i$.
Put $v^+=\big(\mathop{\oplus}\limits_{i\ {\text{even}}}v_i\big)\colon E^+\to E^-$,
and define $w^+$ similarly. In the same way, put $p^\pm=\big(\mathop{\oplus}\limits_{i\ {\text{even/odd}}}p_i\big)\colon E^\pm\to H^\pm$ and define $h^\pm$ similarly, the
class of link to consider is the class of
\begin{equation}\label{linkcohomprime}v^++w^++p^++h^-\colon E^+\oplus H^-\longrightarrow
E^-\oplus H^+\end{equation}
Data as $v$, $w$, $p$, and $h$ can be obtained by choosing hermitian metrics
on the $E^i$ and considering orthocomplements of ${\rm{Im}}v_{i-1}$ in ${\rm{Ker}}v_i$
and of ${\rm{Ker}}v_i$ in $E^i$ which are to be put in bijection with one another
or with $H^i$ (they are all vector subbundles of the $E^i$). The bijectivity of the above map can be proved for any set of $w_i$, $h_i$ and
$p_i$ with required properties, by showing the nullity of the
highest degree component (in the $i$-graduation) of some
element lying in its kernel.
This fact and the fact that the set of $(w,h,p)$-data with
required properties is convex proves the independence of this class of link on
the choices made. Put hermitian metrics on the $E^i$, and consider
the adjoints $v_i^*\colon E^{i+1}\longrightarrow E^i$, and the associated $(v^-)^*=\big(\mathop{\oplus}\limits_{i\ {\text{odd}}}v_i^*\big)\colon E^+\to E^-$. Choose $p_i$ and $h_i$ such that they realise the Hodge isomorphisms
$H^i\cong{\rm{Ker}}(v_i+v_{i-1}^*)$. The link between $E^+
-E^-$ and $H^+-H^-$ is also provided by the following exact sequence
\begin{equation}\label{oscour}0\longrightarrow H^+\overset{h^+}\longrightarrow E^+\overset{v^++(v^-)^*}{-\!\!\!-\!\!\!-\!\!\!-\!\!\!-\!\!\!\longrightarrow}
E^-\overset{p^-}\longrightarrow H^-\longrightarrow0\end{equation}
($p$ being supposed to vanish on ${\rm{Ker}}(v+v^*)^\perp$). This is because the restriction of $v_i\circ v_i^*$ to the image of $v_i$ in $E^i$ is allways isotopic to the identity.
\subsubsection{Flat $K$-theory:}\label{s:flatkth}
For any vector bundle $F$ on $M$ denote by $\Omega^\bullet
(M,F)$ the vector space of smooth differential forms on $M$
with values in $F$.
A connection $\nabla_{\!E}$ on the vector bundle $E$ on $M$ is said
to be flat if its curvature
$\nabla_{\!E}^2\in\Omega^2(M,{\rm{End}}E)$ vanishes. The couple
$(E,\nabla_{\!E})$ is then called a flat vector bundle. Two such flat
bundles $(E,\nabla_{\!E})$ and $(F,\nabla_{\!F})$ are isomorphic if
there exists some smooth vector bundle isomorphism $f\colon
E\to F$ such that $\nabla_{\!E}=f^*\nabla_{\!F}$. Such an $f$ is said to
be a parallel isomorphism.

The group $K^0_{\text{flat}}(M)$ is the quotient of the free
abelian group generated by isomorphism classes of flat vector
bundles, by the following relation:
\begin{equation}\label{suitexacteplate}
(E,\nabla_{\!E})=(E',\nabla_{\!E'})
+(E'',\nabla_{\!E''})\qquad{\text{ if }}\qquad
0\to E'\overset i\longrightarrow E
\overset p\longrightarrow E''\to 0
\end{equation}
is an exact sequence in the category of flat vector bundles,
which means that $\nabla_{\!E}\circ i=i\circ\nabla_{\!E'}$ and
$\nabla_{\!E''}\circ p=p\circ\nabla_{\!E}$ (the morphisms $i$ and $p$ are called parallel in this case).

There is an obvious forgetful map $(E,\nabla_{\!E})\in K^0_{\text{flat}}(M)\longmapsto[E]\in K^0_{\text{top}}(M)$.

For any complex vector bundle $E$ on $M$ endowed with a hermitian metric $h^E$
and a connection $\nabla_{\!E}$, the adjoint transpose connection $\nabla_{\!E}^*$ of $\nabla_{\!E}$
is defined as follows:
\begin{equation}\label{adjointconnection}
h^E({\nabla_{\!E}^*}_{\tt v}\sigma,\theta)={\tt v}.h^E(\sigma,\theta)-h^E(\sigma,{\nabla_{\!E}}_{\tt v}\theta)
\end{equation}
where $\sigma$ and $\theta$ are local sections of $E$, ${\tt v}$ is a tangent
vector, ${\tt v}.f$ is the derivative of the function $f$ along ${\tt v}$,
${\nabla_{\!E}^*}_{\tt v}\sigma$ is the derivative of $\sigma$ along ${\tt v}$ with respect to the connection $\nabla_{\!E}^*$ and accordingly for
${\nabla_{\!E}}_{\tt v}\theta$.
Of course $(\nabla_{\!E}^*)^*=\nabla_{\!E}$, (and $\nabla_{\!E}=\nabla_{\!E}^*$ if and only if
$\nabla_{\!E}$ respects the hermitian metric $h^E$).

The curvatures $\nabla_{\!E}^2$ and $\nabla_{\!E}^{*2}$
are mutually skew adjoint, so that $\nabla_{\!E}^*$ is flat if and only if $\nabla_{\!E}$ is, and this gives rise to a ``conjugation'' involution:
\begin{equation}\label{collier}(E,\nabla_{\!E})\in K^0_{\rm{flat}}(M)\longmapsto (E,\nabla_{\!E}^*)\in
K^0_{\rm{flat}}(M)
\end{equation}
Indeed, for some exact sequence of the form \eqref{suitexacteplate}, its transpose exact sequence
\begin{equation}\label{transposeflat}
0\to E''\overset {p^*}\longrightarrow E
\overset{i^*}\longrightarrow E'\to 0
\end{equation}
turns out to be an exact sequence of flat vector bundles with respect to
transpose adjoint connections on $E''$, $E'$ and $E$. Moreover, if $h^E_1$ and $h^E_2$ are two different hermitian metrics on $E$,
define the global automorphism $g$ of $E$ by
the following formula, valid for any local sections $\sigma$ and $\theta$ of $E$:
\begin{equation}\label{gnabla}
h^E_2(\sigma,\theta)=h^E_1(g(\sigma),\theta)\end{equation}
Call $\nabla_{\!E,1}^*$ and $\nabla_{\!E,2}^*$ the adjoint transpose of $\nabla_{\!E}$ relatively to $h^E_1$ and $h^E_2$ respectively, then
it is easily verified that $\nabla^*_{\!E,2}=g^{-1}\nabla_{\!E,1}^*g$,
so that $(E,\nabla_{\!E,2}^*)$ and $(E,\nabla_{\!E,1}^*)$
are isomorphic (throw $g$) as flat vector bundles, they thus define the same element in $K^0_{\rm{flat}}(M)$.
\subsubsection{Relative $K$-theory:}\label{s:relativeKdefinition}
Consider now on $M$
quintuples $(E,\nabla_{\!E},F,\nabla_{\!F},f)$ where $(E,\nabla_{\!E})$ and
$(F,\nabla_{\!F})$ are
flat vector bundles on $M$, and $f\colon
E\to F$ is a ${\mathcal C}^\infty$ isomorphism.
\begin{definition}\label{definitiondelaKtheorierelative}
$K^0_{\text{rel}}(M)$ is the quotient of the free abelian group
generated by such quintuples modulo the following relations:
\begin{description}
\item[] $(i)$ $(E,\nabla_{\!E},F,\nabla_{\!F},f)=0$ if $f$ is isotopic
to some parallel isomorphism.
\item[] $\begin{aligned}\!(ii)\, \qquad(E,\nabla_{\!E},F,\nabla_{\!F},f)&+(G,\nabla_{\!G},H,\nabla_{\!H},h)=\\&=(E\oplus G,\nabla_{\!E}\oplus\nabla_{\!G},F\oplus H,
\nabla_{\!F}\oplus\nabla_{\!H},f\oplus h)\end{aligned}$
\item[] $(iii)$ $(E,\nabla_{\!E},E'\oplus E'',\nabla_{\!E'}\oplus\nabla_{\!E''}
,s\oplus p)$ vanishes
in $K^0_{\text{rel}}(M)$ if there is a short exact sequence of flat bundles as in \eqref{suitexacteplate} above and if $s\colon E\to E'$
is a ${\mathcal C}^\infty$ bundle map such that $s\circ i$ is the identity
of $E'$.
\end{description}
\end{definition}
\begin{remark}\label{cycle}
Note that $(E,\nabla_{\!E},F,\nabla_{\!F},f)=(E,\nabla_{\!E},F,\nabla_{\!F},g)$ if
$f$ and $g$ are isotopic, that $(E,\nabla_{\!E},F,\nabla_{\!F},f)+(F,\nabla_{\!F},G,\nabla_{\!G},g)=(E,\nabla_{\!E},G,\nabla_{\!G},g\circ f)$, and that $(E,\nabla_{\!E},F,\nabla_{\!F},f)=
(E',\nabla_{\!E'},F',\nabla_{\!F'},f')+
(E'',\nabla_{\!E''},F'',\nabla_{\!F''},f'')\ $ if
\[\begin{CD}
0@>>>E'@>>>E@>>>E''@>>>0\\
@.@V{f'}VV@VVfV@VV{f''}V\\
0@>>>F'@>>>F@>>>F''@>>>0
\end{CD}\]
is a commutative diagram whose lines are short exact sequences in
the category of flat vector bundles (on $M$).

In fact the first one and the third one of these three relations are together equivalent to $(i)$,
$(ii)$ 
and $(iii)$ so that they can be used to provide an alternative
definition of $K^0_{\text{rel}}(M)$ (see \cite{BerthomieuPNCDIRKT} \S2.1 for details).

Now if $(E,\nabla_{\!E})$, $(F,\nabla_{\!F})$, $(G,\nabla_{\!G})$ and $(H,\nabla_{\!H})$ are flat vector
bundles on $M$, then any link $({\mathbb C}^N,\ell)$ between $E-F$ and $G-H$ defines
an element of $K^0_{\text{rel}}(M)$ as
$(E\oplus H\oplus{\mathbb C}^N,\nabla_{\!E}\oplus\nabla_{\!H}\oplus d_{{\mathbb C}^N},F\oplus G\oplus{\mathbb C}^N,\nabla_{\!F}\oplus\nabla_{\!G}\oplus d_{{\mathbb C}^N},\ell)$, where $d_{{\mathbb C}^N}$
is the canonical flat connections on the trivial vector bundle ${\mathbb C}^N$. Two equivalent links will provide equivalent elements of $K^0_{\text{rel}}(M)$, (and any link
between $E-F$ and $G-H$ is equivalent to one of the form $({\mathbb C}^N,\ell)$). This unambiguously defined element of $K^0_{\rm{rel}}(M)$ will be denoted by $(E\oplus H,\nabla_{\!E}\oplus\nabla_{\!H},F\oplus G,\nabla_{\!F}\oplus\nabla_{\!G},[\ell])$.

Finally, for exact sequences in the category of flat bundles as
\begin{equation}\label{E:sixdiagram}
\begin{aligned}0\to E\longrightarrow G\longrightarrow H&\longrightarrow F\to 0\\
{\text{or }}\qquad\begin{CD}E@>>>G\\@AAA@VVV\\
F@<<<H\end{CD}\qquad{\text{ or }}\qquad&
\begin{CD}A@>>>B@>>>C\\@AAA@.@VVV\\F@<<<E@<<<D\end{CD}\qquad\end{aligned}
\end{equation}
if $[\ell]$ and $[\ell']$ are the equivalence classes of links between $E-F$ and $G-H$
(in either of the two first cases) or between $B-E$ and
$(A\oplus C)-(D\oplus F)$ respectively, which are defined by the corresponding exact sequences (as in \eqref{precedente} and \eqref{complexestousrondsavant} before),
then the following associated elements
\begin{align*}
(E\oplus H,\nabla_{\!E}\oplus\nabla_{\!H},G&\oplus F,\nabla_{\!G}\oplus\nabla_{\!F},[\ell])\\
(B\oplus D\oplus F,\nabla_{\!B}\oplus\nabla_{\!D}\oplus\nabla_{\!F},
A&\oplus C\oplus E,
\nabla_{\!A}\oplus\nabla_{\!C}\oplus\nabla_{\!E},[\ell'])
\end{align*}
vanish in $K_{\text{rel}}^0(M)$, as can be proved by decomposing these exact sequences in short ones. Note that the fact that the morphisms are parallel forces them to be of constant
rank in contrast with the situation of \eqref{complexestousrondsavant}.
\end{remark}
One has the following morphisms:
\[\begin{aligned}
\big[g:{\mathbb C}^N\overset\sim\longrightarrow{\mathbb C}^N\big]\in K_1^{\text{top}}(M)&\longmapsto({\mathbb C}^N,d_{{\mathbb C}^N},
{\mathbb C}^N,d_{{\mathbb C}^N},g)\in K_0^{\text{rel}}(M)\\
(E,\nabla_{\!E},F,\nabla_{\!F},f)\in K_0^{\text{rel}}(M)&\overset\partial\longmapsto (F,\nabla_{\!F})-(E,\nabla_{\!E})\in K_0^{\text{flat}}(M)
\end{aligned}\]
which enter in the following exact sequence (see \cite {BerthomieuPNCDIRKT} \S2.2 for a proof):
\begin{equation}\label{suitexmoi}
K^1_{\text{top}}(M)\longrightarrow K^0_{\text{rel}}(M)\overset\partial\longrightarrow K^0_{\text{flat}}(M)\longrightarrow K^0_{\text{top}}(M)
\end{equation}

There is the same conjugation involution in $K^0_{\rm{rel}}(M)$ as in $K^0_{\rm{flat}}(M)$:
\begin{equation}\label{congrel}(E,\nabla_{\!E},F,\nabla_{\!F},f)\longmapsto (E,\nabla_{\!E}^*,F,\nabla_{\!F}^*,f)\end{equation}
where $\nabla_{\!E}^*$ and $\nabla_{\!F}^*$ are adjoint transpose of $\nabla_{\!E}$ and $\nabla_{\!F}$ with respect to any hermitian metrics on $E$ and $F$. This is due firstly to the independence on these metrics of the image as element of $K^0_{\rm{rel}}(M)$, (which can be checked in the same way as was made in the case of the
conjugation in $K^0_{\rm{flat}}$) and secondly to the equivalence of relation $(iii)$
in definition \ref{definitiondelaKtheorierelative} with the following:
if there is a short exact sequence of flat bundles as in \eqref{suitexacteplate} above and if $j\colon E''\to E$
is a ${\mathcal C}^\infty$ bundle map such that $p\circ j$ is the identity
of $E''$, then $(E'\oplus E'',\nabla_{\!E'}\oplus\nabla_{\!E''}
,E,\nabla_{\!E},i+j)$ vanishes
in $K^0_{\text{rel}}(M)$ (in fact $i+j$ is isotopic to $(s\oplus p)^{-1}$, and
$s^*$ in \eqref{transposeflat} is in the same position as $j$ here).

Elements of $K^0_{\rm{rel}}(M)$ of the form $(E,\nabla_{\!E},E,\nabla_{\!E}^*,
{\rm{Id}}_E)$ are purely imaginary with respect to this conjugation; reciprocally, the subgroup of $K^0_{\rm{rel}}(M)$ generated by such elements
is equal to, or of index 2 in, the purely imaginary part of $K^0_{\rm{rel}}(M)$.
This is because (see the beginning of remark \ref{cycle})
\begin{equation}\begin{aligned}
(E,\nabla_{\!E},F,\nabla_{\!F},f)-&(E,\nabla_{\!E}^*,F,\nabla_{\!F}^*,f)=
\\&=(F,\nabla_{\!F},F,\nabla_{\!F}^*,{\rm{Id}}_F)-(E,\nabla_{\!E},E,\nabla_{\!E}^*,{\rm{Id}}_E)\end{aligned}\end{equation}  
\subsection{Chern-Simons transgression:}
\subsubsection{The case of the Chern character:}
Consider some smooth vector bundle $E$ with connection
$\nabla$ on $M$. Chern-Weil theory
associates to such data the following complex differential form on $M$
\begin{equation}\label{premierChern}{\rm{ch}}(\nabla)={\rm{Tr}}\exp\left(-\frac1{2\pi i}
\nabla^2\right)=\phi{\rm{Tr}}\exp(-\nabla^2)\end{equation}
where $\phi$ is the operator on even degree differential forms which divides $2k$-degree forms by $(2\pi i)^k$.
This form is closed, its de Rham cohomology class is
independent of $\nabla$ and equals the image of the Chern
character of $E$ in $H^{\text{even}}(M,{\mathbb C})$.

Note that if $\nabla_{\!E}$ respects some hermitian metric $h^E$
on $E$, then its curvature $\nabla_{\!E}^2$ is skew adjoint, so that
${\rm{ch}}(\nabla_{\!E})$ is a real differential form.

If now $\nabla_{\!E,0}$ and $\nabla_{\!E,1}$ are two connections on
the complex vector bundle $E$,
Chern-Simons theory associates to them an odd differential form
$\widetilde{\rm{ch}}(\nabla_{\!E,0},\nabla_{\!E,1})$ defined modulo
coboundaries which verifies
\[d\widetilde{\rm{ch}}(\nabla_{\!E,0},\nabla_{\!E,1})={\rm{ch}}(\nabla_{\!E,1})-
{\rm{ch}}(\nabla_{\!E,0})\]
This form is functorial by pull-backs, additive for direct sums of vector bundles,
it verifies modulo exact forms:
\begin{equation}\label{E:cocycleChS}\widetilde{{\rm{ch}}}(\nabla_{\!E,0},\nabla_{\!E,2})=\widetilde{\rm{ch}}(\nabla_{\!E,0},\nabla_{\!E,1})
+\widetilde{\rm{ch}}(\nabla_{\!E,1},\nabla_{\!E,2})\end{equation}
and it is locally gauge invariant, which means that $\widetilde{\rm{ch}}
(\nabla,g^*\nabla)$ is exact if $g$ is a global ${\mathcal C}^\infty$
automorphism of $E$ isotopic to the identity.

The form $\widetilde{\rm{ch}}(\nabla_{\!E,0},\nabla_{\!E,1})$ is calculated as follows:
consider $p_1\colon M\times[0,1]
\to M$ (the projection on the first factor) and the bundle $p_1^*E$ on $M\times[0,1]$,
choose any connexion
$\widetilde\nabla_{\!E}$ on it such that $\widetilde\nabla_{\!E}\vert_{M\times\{i\}}=\nabla_{\!E,i}$
for $i=0$ and $1$, denote by $\nabla_{\!E,t}=\widetilde\nabla_{\!E}\vert_{M\times\{t\}}$ the fibral restrictions of $\widetilde\nabla$ for all $t\in[0,1]$, and extend $\phi$ to odd forms by deciding that $\phi$ divides $(2k-1)$-degree forms by $(2\pi i)^k$, then
\begin{equation}\label{explicitransgression}
\begin{aligned}
\widetilde{\rm{ch}}(\nabla_{\!E,0},\nabla_{\!E,1})&=\int_{[0,1]}{\rm{ch}}(\widetilde\nabla_{\!E})=-\int_0^1\phi{\rm{Tr}}\left(\frac{\partial\nabla_{\!E,t}}{\partial t}\exp\left(-\nabla_{\!E,t}^2\right)\right)dt\\
&=-\frac1{2\pi i}\int_0^1{\rm{Tr}}\left(\frac{\partial\nabla_{\!E,t}}{\partial t}\exp\left(-\frac1{2\pi i}\nabla_{\!E,t}^2\right)\right)dt
\end{aligned}\end{equation}
The such defined form depends on $\widetilde\nabla_{\!E}$ only throw addition of an exact form.
\subsubsection{The case of the Euler class:}\label{Euler}

In the case of a real vector bundle $F_{\mathbb R}$ endowed with an euclidean metric and
an unitary connection $\nabla_{\!F_{\mathbb R}}$, the curvature $\nabla_{\!F_{\mathbb R}}^2$
is a two-form with values in antisymmetric endomorphisms of
$F_{\mathbb R}$. Define $e(\nabla_{\!F_{\mathbb R}})$ to be zero if $F_{\mathbb R}$ is of odd rank (as
real vector bundle) and to be the P\!faffian of $\frac1{2\pi}
\nabla_{\!F_{\mathbb R}}^2$ if $F_{\mathbb R}$ is of even rank. One obtains a closed real
differential form whose degree equals the rank of $F_{\mathbb R}$ and whose de Rham cohomology class is independent on $\nabla_{\!F_{\mathbb R}}$ (and on the euclidean metric on $F_{\mathbb R}$) and coincides with the image of the Euler class of $F_{\mathbb R}$ in $H^\bullet(M,{\mathbb C})$.

If now $\nabla_{\!F_{\mathbb R},0}$ and $\nabla_{\!F_{\mathbb R},1}$ are two unitary connections 
(with respect to not necessarily same metrics) on $F_{\mathbb R}$, then take some
unitary connection $\widetilde\nabla_{\!F_{\mathbb R}}$ on $p_1^*F_{\mathbb R}$ on $M\times
[0,1]$ such that $\widetilde\nabla_{\!F_{\mathbb R}}\vert_{M\times i}=\nabla_{\!F_{\mathbb R},i
}$ for $i=0$ and $1$ (of course the same restriction property for the metric on $p_1^*F_{\mathbb R}$ is needed), then the form
\begin{equation}\widetilde e(\nabla_{\!F_{\mathbb R},0},\nabla_{\!F_{\mathbb R},1})=\int_{[0,1]}e(
\widetilde\nabla_{\!F_{\mathbb R}})\end{equation}
depends on $\widetilde\nabla_{\!F_{\mathbb R}}$ throw addition of an exact form and verifies
\begin{equation}d\widetilde e(\nabla_{\!F_{\mathbb R},0},\nabla_{\!F_{\mathbb R},1})=e(\nabla_{\!F_{\mathbb R},1})
-e(\nabla_{\!F_{\mathbb R},0})\end{equation}
Its class modulo exact forms is functorial by pullback, locally gauge invariant, and verifies a formula analogous to \eqref{E:cocycleChS}.

Now making the product of $e(\widetilde\nabla_{\!F_{\mathbb R}})$ and ${\rm{ch}}
(\widetilde\nabla_{\!E})$ yields the following equality modulo exact forms:
\begin{equation}\label{E:wedge}\begin{aligned}
\int_{[0,1]}e(\widetilde\nabla_{\!F_{\mathbb R}})\wedge {\rm{ch}}(\widetilde\nabla_{\!E})&=\widetilde{e}(\nabla_{\!F_{\mathbb R},0},
\nabla_{\!F_{\mathbb R},1})\wedge{\rm{ch}} (\nabla_{\!E,0})+e(\nabla_{\!F_{\mathbb R},1})\wedge\widetilde{\rm{ch}}
(\nabla_{\!E,0},\nabla_{\!E,1})\\
&=e(\nabla_{\!F_{\mathbb R},0})\wedge\widetilde{\rm{ch}}(\nabla_{\!E,0},
\nabla_{\!E,1})
+\widetilde{e}(\nabla_{\!F_{\mathbb R},0},\nabla_{\!F_{\mathbb R},1})\wedge{\rm{ch}}(\nabla_{\!E,1})
\end{aligned}\end{equation}
\subsubsection{Chern-Simons transgression and hermitian metrics:}\label{adjcon}
For a complex vector bundle $E$ on $M$ endowed with an hermitian
metric $h^E$ and a connection $\nabla_{\!E}$, remember the definition of the 	adjoint transpose connection $\nabla_{\!E}^*$ from \eqref{adjointconnection}.
The fact that the curvatures of $\nabla_{\!E}$ and $\nabla_{\!E}^*$
are mutually skew adjoint has the following consequence:
\begin{equation}\label{conjugateChernchar}{\rm{ch}}(\nabla_{\!E}^*)=\overline{{\rm{ch}}(\nabla_{\!E})}\end{equation}
\begin{lemma}\label{indepmetrique}
The class of $\widetilde{\rm{ch}}(\nabla_{\!E}^*,\nabla_{\!E})$ in
$\Omega(M,{\mathbb C})/d\Omega(M,{\mathbb C})$ is independent
of the metric $h^E$.
\end{lemma}
\begin{proof}
This is a simple consequence of the local gauge invariance
of the Chern-Simons forms modulo exact forms. For two metrics
$h^E_0$ and $h^E_1$, choose a path $h^E_t$ of hermitian metrics connecting them, and define for any $t$ the global automorphism $g_t$ of $E$ by the following formula, valid for any sections $\sigma$ and $\theta$ of $E$: 
\[h^E_t(\sigma,\theta)=h^E_0(g_t\sigma,\theta)\]
Then the adjoint transpose of $\nabla_{\!E}$ with respect to
$h^E_t$ equals $g_t^{-1}\nabla_{\!E}^*g_t$ so that
it is a gauge transform of $\nabla_{\!E}^*$ with gauge $g_t$, which
is isotopic to $g_0={\rm{Id}}$
\end{proof}
If $\nabla_{\!E,0}$ and $\nabla_{\!E,1}$ are connections on $E$,
then using \eqref{conjugateChernchar} above and considering the adjoint
transpose connection of $\widetilde\nabla_{\!E}$ in formula
\eqref{explicitransgression} yields the following relation modulo exact forms:
\begin{equation}\label{metricandconnection}
\widetilde{\rm{ch}}(\nabla_{\!E,0}^*,\nabla_{\!E,1}^*)=
\overline{\widetilde{\rm{ch}}(\nabla_{\!E,0},\nabla_{\!E,1})}
\end{equation}
where $\nabla^*_{\!E,0}$ and $\nabla^*_{\!E,1}$ are adjoint transpose of $\nabla_{\!E,0}$ and $\nabla_{\!E,1}$ with respect to possibly different hermitian metrics on $E$. In particular
$\widetilde{\rm{ch}}(\nabla_{\!E,0},\nabla_{\!E,1})$ is a
real form if $\nabla_{\!E,0}$ and $\nabla_{\!E,1}$ are
unitary connection (with respect to possibly different metrics on $E$), and in general
\[\widetilde{\rm{ch}}(\nabla_{\!E}^*,\nabla_{\!E})=
-\widetilde{\rm{ch}}(\nabla_{\!E},\nabla_{\!E}^*)=-\overline{
\widetilde{\rm{ch}}(\nabla_{\!E}^*,\nabla_{\!E})}\]
is a purely imaginary form (modulo exact forms).

Consider the connection
$\nabla_{\!E}^u=\frac12(\nabla_{\!E}+\nabla^*_{\!E})$; it respects $h^E$, and
\begin{equation}\label{nonflatimaginary}
\widetilde{\rm{ch}}(\nabla_{\!E}^*,\nabla_{\!E})=
\widetilde{\rm{ch}}(\nabla_{\!E}^*,\nabla_{\!E}^u)+
\widetilde{\rm{ch}}(\nabla_{\!E}^u,\nabla_{\!E})=
2i{\mathfrak{Im}}\big(\widetilde{\rm{ch}}
(\nabla_{\!E}^u,\nabla_{\!E})\big)
\end{equation}
(Note that if $\nabla_{\!E}$ is flat, then $\nabla_{\!E}^*$ is flat too, but
$\nabla_{\!E}^u$ need not be flat, so that working with flat connections forces one to work with nonmetric connections).

Finally, if $\nabla_{\!E,0}$ and $\nabla_{\!E,1}$ are connections on $E$, then the cocycle condition \eqref{E:cocycleChS} produces the following
relation modulo exact forms
\begin{equation}\label{boborelrel}
\begin{aligned}\widetilde{\rm{ch}}(\nabla_{\!E,1}^*,\nabla_{\!E,1})
-\widetilde{\rm{ch}}(\nabla_{\!E,0}^*,\nabla_{\!E,0})
&=\widetilde{\rm{ch}}(\nabla_{\!E,0},\nabla_{\!E,1})
-\widetilde{\rm{ch}}(\nabla_{\!E,0}^*,\nabla_{\!E,1}^*)
\\&=2i\, {\mathfrak{Im}}\, \widetilde{\rm{ch}}(\nabla_{\!E,0},
\nabla_{\!E,1})
\end{aligned}\end{equation}
\begin{remark}\label{Lottremark}
If $\nabla_{\!E}$ is flat, then it is proved in \cite{BismutLott}, proof of proposition 1.14 that in fact $\widetilde{\rm{ch}}
(\nabla_{\!E},\nabla_{\!E}^u)$ is purely imaginary so that
\eqref{nonflatimaginary} holds without $i\, {\mathfrak{Im}}$. Moreover 
\[\widetilde{\rm{ch}}(\nabla_{\!E}^*,\nabla_{\!E})=\frac1\pi\sum_{j=0}^\infty\frac{2^{2j}j!}{(2j+1)!}c_{2j+1}(E)\]
is the degree decomposition of $\widetilde{\rm{ch}}(\nabla_{\!E}^*,\nabla_{\!E})$, where the $c_k(E)$ are the classes considered by Bismut and Lott (see \cite{BismutLott} formulae (0.2) and (1.34)). These are exactly the imaginary part of the Chern-Cheeger-Simons classes of flat complex vector bundles (\cite{CheegerSimons} and \cite{BismutLott} Proposition 1.14).
\end{remark}
\begin{lemma}\label{4k+1ou3}
If $E_{\mathbb R}$, is a real vector bundle on
$M$ with connection $\nabla_{\!E_{\mathbb R}}$, and if $E$ is its complexification with associated connection $\nabla$, then, up to exact forms, $\widetilde{\rm{ch}}(\nabla^*,\nabla)$ vanishes in degrees $4k+3$. 
\end{lemma}
\begin{proof} Suppose that $E$ is endowed with a hermitian form which is the
complexification of a real scalar product on $E_{\mathbb R}$, and use
the path of connections $\nabla_{\!t}=(1-t)\nabla^*+t\nabla$, then
the lemma follows from formula \eqref{explicitransgression} by counting the $i$,
and from lemma \ref{indepmetrique}.
\end{proof}
\subsubsection{Superconnections:}\label{SuperK}
Consider a complex vector space $V$ which decomposes as 
$V=V^+\oplus V^-$, with a ${\mathbb Z}_2$-graduation operator
$\tau\vert_{V^\pm}=\pm{\rm{Id}}\vert_{V^\pm}$.
The supertrace of $a\in{\rm{End}}V$ is defined by
${\rm{Tr}}_sa={\rm{Tr}}(\tau\circ a)$, (this is the trace on $V^+$
minus the trace on $V^-$).
${\rm{End}}V$ is also ${\mathbb Z}_2$-graded (even endomorphisms respect
both parts $V^+$ and $V^-$ and odd ones exchange them). The supercommutator in ${\rm{End}}V$
is defined for pure degree objects as
\[[a,b]=ab-(-1)^{{\rm{deg}}a\, {\rm{deg}}b}ba\]
and bilinearly extended to ${\rm{End}}V$. This is such that the supertrace vanishes on supercommutators. In particular, any endomorphism which commutes
(or supercommutes, because odd endomorphisms are traceless)
with some odd invertible endomorphism (if there exists some...) has vanishing supertrace.

Suppose now that $E$ is a ${\mathbb Z}_2$ graded vector bundle on $M$,
that is $E=E^+\oplus E^-$ where $E^+$ and $E^-$ are complex vector bundles themselve.
The supertrace is defined as above and extends naturally on ${\rm{End}}E$-valued
differential forms, with values in ordinary differential forms.
${\rm{End}}E$-valued differential forms inherit a global
${\mathbb Z}_2$-graduation, ordinary differential forms being ${\mathbb Z}_2$-graduated by the parity of their degree. They act on $E$-valued differential forms
and multiply in the following way
\begin{equation}(\alpha\widehat\otimes a)(\beta\widehat\otimes b)=(-1)^{{\rm{deg}}a\, {\rm{deg}}
\beta}(\alpha\wedge\beta)\otimes(ab)\end{equation}
where $\alpha\widehat\otimes a$ and $\beta\widehat\otimes b$ are decomposed
tensors in the graded tensor product of differential forms with either
${\rm{End}}E$ or $E$. The supercommutator of ${\rm{End}}E$-valued
differential forms is defined in the same way as above but by considering
the global graduation. With this convention, the supertrace allways
vanishes on supercommutators, and ${\rm{End}}E$-valued
differential forms which supercommute with some odd invertible endomorphism
of $E$ have vanishing supertrace.

A superconnection $A$ on $E$ is the sum of a connection $\nabla$
which respects the
decomposition of $E$ and of a globally odd ${\rm{End}}E$-valued differential
form $\omega$. Its curvature is its square
$A^2=(\nabla+\omega)^2=\nabla^2+[\nabla,\omega]+\omega^2$, a global even ${\rm{End}}E$-valued
differential form ($A^2$ is not a differential operator).

Following \eqref{premierChern}, denote by ${\text{ch}}(A)=\phi{\text{Tr}}_{s}
\exp-A^2$ the Chern-Weil form representing the Chern character of
any superconnection $A$. It is an even
degree differential form on $M$. The space of superconnections on $E$
is convex (and of course contains ordinary connections) so that the
preceding Chern-Weil and Chern-Simons theory also works for superconnections
(especially formula \eqref{explicitransgression}).
Thus ${\rm{ch}}(A)$ is closed and its cohomology class is the
same as the Chern character of $E$ in complex cohomology (this means ${\rm{ch}}(E^+)-{\rm{ch}}(E^-)$ because of the ${\mathbb Z}_2$-graduation).

A hermitian metric on $E=E^+\oplus E^-$ will be supposed to
make this decomposition orthogonal. Let $\beta$ be a differential form
and $a\in{\rm{End}}E$, the adjoint of $a$ will be denoted by
$a^*$. For ${\rm{End}}E$-valued differential forms, there are
two notions of adjunction: the ordinary adjoint of $\beta\widehat\otimes a$
will be $\overline\beta\widehat\otimes a^*$, while its
special adjoint will be $(-1)^{\frac{{\rm{deg}}\beta({\rm{deg}}\beta-1)}2+{\rm{deg}}\beta{\rm{deg}}a}
\overline\beta\widehat\otimes a^*$, following
the convention implicitely used in \cite{BismutLott} \S I(c) and (d). If $\omega_1$ and $\omega_2$ are any (multidegree)
${\rm{End}}E$-valued differential forms, denote by $\omega_1^S$
and $\omega_2^S$ their special adjoints, then $(\omega_1\omega_2)^S=\omega_2^S\omega_1^S$.

As usual, the trace (or supertrace) of an autoadjoint ${\rm{End}}E$-valued
differential form is a real form and the trace (or supertrace) of a
skewadjoint one is a purely imaginary form, while the trace (or supertrace)
of a special autoadjoint one is real in degrees $4k$ and $4k-1$ and purely
imaginary in degrees $4k+1$ and $4k+2$.
Consequently, $\phi{\rm{Tr}}_s(\omega)$ is real if $\omega$ is a special autoadjoint (multidegree) ${\rm{End}}E$-valued differential form (and in general, $\phi{\rm{Tr}}_s(\omega)=\overline{\phi{\rm{Tr}}_s(\omega^S)}$).

Let $\omega$ be some globally odd ${\rm{End}}E$-valued differential form,
and $A=\nabla+\omega$ a superconnection on $E$. Then the adjoint transpose of $A$ is
\begin{equation}A^S=\nabla^*+\omega^S\end{equation}
Thus $\frac12(A+A^S)$ is the sum of $\nabla^u$ (which respects the hermitian metric of $E$)
and of some special autoadjoint ${\rm{End}}E$-valued differential form of globally odd degree. The fact that $[\nabla^*,\omega^S]=[\nabla,\omega]^S$
and $(\nabla^*)^2=-(\nabla^2)^*=(\nabla^2)^S$ has the consequence that
$(A^S)^2=(A^2)^S$ so that ${\rm{ch}}(A^S)=\overline{{\rm{ch}}(A)}$.
In particular, ${\rm{ch}}\big(\frac12(A+A^S)\big)$ is a real form.
Finally, lemma \ref{indepmetrique} and formulae \eqref{metricandconnection} and \eqref{nonflatimaginary} also hold in the context of superconnections.
\subsubsection{Chern-Simons transgression, exact sequences and links:}
Let $E$, $F$, $G$ and $H$ be vector bundles on $M$ with connections
$\nabla_{\!E}$, $\nabla_{\!F}$, $\nabla_{\!G}$ and $\nabla_{\!H}$, and
suppose there exists some link $(K,\ell)$ between $E-F$ and $G-H$.
One can then associate to $(K,\ell)$ the
differential form (defined modulo exact forms)
\begin{equation}\label{E:translink}
\widetilde{\text{ch}}([\ell])=\widetilde{\text{ch}}\big(\nabla_{\!E}\oplus
\nabla_{\!H}\oplus \nabla_{\!K},\ell^*[\nabla_{\!F}\oplus\nabla_{\!G}\oplus\nabla_{\!K}]\big)
\end{equation}
for some
connection $\nabla_{\!K}$ on $K$. It is easily checked that
the class of this form modulo exact forms does not depend on
the choice of $\nabla_{\!K}$ and is not modified by changing $(K,\ell)$ by
an equivalent link. It is possible to choose a unitary $\nabla_{\!K}$, so that $\widetilde{\rm{ch}}([\ell])$ is a real form (modulo exact forms)
if it happens that $\nabla_{\!E}$, $\nabla_{\!F}$, $\nabla_{\!G}$ and $\nabla_{\!H}$ are all unitary connections.

For the composition of two links $\ell$ and $\ell'$, and any connections on the
considered bundles one obtains (modulo exact forms):
\begin{equation}\label{E:gof}
\widetilde{\text{ch}}([\ell'\circ \ell])=\widetilde{\text{ch}}([\ell])
+\widetilde{\text{ch}}([\ell'])
\end{equation}

Now consider a short exact sequence of vector bundles on $M$
\[0\to E'\overset i\longrightarrow E\overset p\longrightarrow E''\to 0\]
and connections $\nabla_{\!E'}$, $\nabla_{\!E}$ and $\nabla_{\!E''}$ on $E'$, $E$ and $E''$ (not necessarily flat) for which $i$ and $p$ are parallel. Let $s\colon E\to E'$ be any bundle morphism such that $s\circ i={\rm{Id}}_{E'}$ then
\begin{lemma}\label{L:exactseqtransgr} The form $\widetilde{\rm{ch}}\big(\nabla_{\!E},
(s\oplus p)^*(\nabla_{\!E'}\oplus\nabla_{\!E''})\big)$ vanishes.
\end{lemma}
\begin{proof}$s$ provides an isomorphism $s\oplus p\colon E{\overset\sim\longrightarrow} E'\oplus
E''$, and the fact that $i$ and $p$ are parallel means that
with respect to this decomposition, the connexions $\nabla_{\!E}$
and $(s\oplus p)^*(\nabla_{\!E'}\oplus\nabla_{\!E''})$ differ from a
one-form $\omega$ with values in ${\rm{Hom}}(E'',E')$.

Consider the path of connexions $\nabla_{\!t}=\nabla_{\!E}-t\omega$.
Of course $\nabla_{\!t}^2$ is upper triangular with respect to the
decomposition $E\cong E'\oplus E''$ and $\omega$ also is, but with
vanishing diagonal terms. The trace thus vanish in formula \eqref{explicitransgression} applied to our situation, and this proves the lemma.
\end{proof}
One consequence of this lemma is the nullity of the classes $\widetilde{\rm{ch}}$
of the links produced by exact sequences as in
formula \eqref{precedente} or \eqref{complexestousrondsavant} if it happens that the bundles
are endowed with such connections that the morphisms be parallel
(this and the fact that kernels and images of the morphisms are of constant rank is true from the parallelism of the morphisms, the connections
being flat as in \eqref{E:sixdiagram} or not). In the same way,
for a complex of vector bundles as described in \eqref{E:linkcohom} (from which the notations are taken), the class
$\widetilde{\rm{ch}}$ of the obtained link between $E^+-E^-$ and $H^+-H^-$  vanishes if the connections on the $E^i$ are such
that the $v_i$ are parallel, and if the connections on the $H^i$ are those
obtained from the quotient $H^i\cong{\rm{Ker}}v_i/{\rm{Im}}v_{i-1}$
(the parallelism of the $v_i$ allows to define such quotient connections).
\subsubsection{``Nadel's'' classes on $K_0^{\text{rel}}(M)$:}
\begin{definition}\label{Nadel}
\[{\mathcal N}_{\rm{ch}}(E,\nabla_{\!E},F,\nabla_{\!F},f)=
\left[\widetilde{\rm{ch}}(\nabla_{\!E},
f^*\nabla_{\!F})\right]\in H^{\text{odd}}(M,{\mathbb C})\]
(of course $\widetilde{\rm{ch}}(\nabla_{\!E},
f^*\nabla_{\!F})$ is closed since ${\rm{ch}}(\nabla_{\!E})$ and ${\rm{ch}}(f^*\nabla_{\!F})$
both equal ${\rm{rk}}E$).
\end{definition}
Arguments as in \cite{BerthomieuPNCDIRKT} Theorem 3.5 and its Corollary allow to prove the
\begin{theorem*}
${\mathcal N}_{\rm{ch}}$ induces a group morphism from $K_0^{\text{rel}}(M)$ to $H^{\text{odd}}
(M,{\mathbb C})$.
\end{theorem*}
Arguments as in \cite{BerthomieuPNCDIRKT} \S 5.1 and 5.2 or \cite{Nadel} allow to prove the following facts:
\begin{description}
\item[]
Let $\Phi$ multiply $2k$ and $(2k-1)$-degree forms by $k!$, then
$\Phi{\rm{ch}}$ is the Chern character without denominators.
The nonintegrity of $\Phi{\mathcal N}_{\rm{ch}}(E,\nabla_{\!E},F,\nabla_{\!F},f)$ detects the fact that $(E,\nabla_{\!E})\neq(F,\nabla_{\!F})\in K_0^{\text{flat}}
(M)$.
\item[]
The nonintegrity of the degree$\geq3$ components of $\Phi{\mathcal N}_{\rm{ch}}(E,\nabla_{\!E},F,\nabla_{\!F},f)$ detects the fact that $(F,\nabla_{\!F})$ cannot be obtained from $(E,\nabla_{\!E})$
throw a deformation of flat bundles, where a deformation of flat bundles on $M$ is a ${\mathcal C}^\infty$
vector bundle $\widetilde E$ on $M\times[0,1]$ with a connection $\widetilde\nabla$ whose restriction
to $E_t=\widetilde E\vert_{M\times\{t\}}$ is flat for any point $t\in [0,1]$ and such that $(E_0,\widetilde\nabla\vert_{M\times\{0\}})\cong(E,\nabla_{\!E})$
and $(E_1,\widetilde\nabla\vert_{M\times\{1\}})\cong(F,\nabla_{\!F})$
\item[]
The nonnullity of the degree$\geq3$ components of ${\mathcal N}_{\rm{ch}}(E,\nabla_{\!E},F,\nabla_{\!F},f)$ detects the fact that $(F,\nabla_{\!F})$ cannot be obtained from $(E,\nabla_{\!E})$
throw a deformation of flat bundles, for which the parallel transport along $[0,1]$ would be isotopic to $f$.
\item[]If $(F,\nabla_{\!F})$ can be obtained from
$(E,\nabla_{\!E})$ by deformation of flat bundles, then
the degree 1 component of ${\mathcal N}_{\rm{ch}}$
modulo integral cohomology
detects the variation of the determinant line. 
\end{description}
The third statement is known as the rigidity of higher classes of flat bundles.
\begin{remark*}
Let $(E,\nabla_{\!E},F,\nabla_{\!F},f)\in K^0_{\text{flat}}(M)$, and consider
the following one-form with values in endomorphisms
of $E$: $\omega=f^*\nabla_{\!F}-\nabla_{\!E}$ . Then it can be proved as in \cite{BerthomieuPNCDIRKT} lemma 4.3
that in fact
\[\widetilde{\rm{ch}}(\nabla_{\!E}, f^*\nabla_{\!F})=-\sum_{r=1}^{\left[\frac{{\rm{dim}}M}2\right]}\left(\frac1{2\pi i}\right)^r\frac{(r-1)!}{(2r-1)!}{\rm{Tr}}(\omega^{2r-1})\]
so that ${\mathcal N}_{\rm{ch}}(E,\nabla_{\!E},F,\nabla_{\!F},f)$ can be computed in the same way as the classes studied in \cite{Nadel} and \cite{BismutLott}.
(This is of course a particular property of flat connections and cannot be
generalised to any connections).

Another remark is that purely imaginary elements of $K^0_{\rm{rel}}(M)$
have purely imaginary ${\mathcal N}_{\rm{ch}}$ classes (see \eqref{congrel} and \eqref{metricandconnection}).
\end{remark*}
\subsection{Free multiplicative $K$-theory:}
\subsubsection{A $K$-theory group with Chern-Weil Chern character:}
Consider some triple $(E,\nabla_{\!E},\alpha)$ where $E$ is a
smooth complex vector bundle on $M$, $\nabla_{\!E}$ a connexion
on $E$ and $\alpha$ is an odd degree differential form defined
modulo exact forms. Two such objects
$(E_1,\nabla_{\!E_1},\alpha_1)$ and $(E_2,\nabla_{\!E_2},\alpha_2)$
will be equivalent if there is some smooth vector bundle
isomorphism $f\colon E_1\to E_2$ such that
\begin{equation}\label{Kahatrelation}
\alpha_2=\alpha_1+\widetilde{\rm{ch}}
(\nabla_{\!E_1},f^*\nabla_{\!E_2})
\end{equation}
The following free multiplicative $K$-group is considered by Bunke and Schick \cite{BunkeSchick}:
\begin{definition}\label{D:deftransgrKth}
The group of free multiplicative $K$-theory $\widehat K_{\rm{ch}}(M)$
is the quotient of the free abelian group generated by triples
$(E,\nabla_{\!E},\alpha)$ as above modulo the preceding equivalence
relation and modulo direct sum (of the vector bundle, with
direct sum connection and sum of the differential forms).

The Chern character on $\widehat K_{\rm{ch}}(M)$ is the map
\[\widehat{\rm{ch}}\colon (E,\nabla_{\!E},\alpha)\in\widehat
K_{\rm{ch}}(M)\, \longmapsto\, {\rm{ch}}(\nabla_{\!E})
-d\alpha\in\Omega^{\rm{even}}(M)/d\Omega^{\rm{odd}}(M)\]

The kernel of the above Chern character will be denoted $K^{-1}_{{\mathbb C}/{\mathbb Z}}(M)$ following \cite{Lott}
Definition 3, and the preimage of ${\mathbb Z}$ by this Chern
character is the (Karoubi) mutiplicative $K$-theory group
$MK^0(M)$ we will consider here (see \cite{KaroubiHCKT} \S 7.5
and \cite{KaroubiCCFFHA} EXEMPLE 3).
\end{definition}
Of course, $MK^0(M)\cong{\mathbb Z}\oplus K^{-1}_{{\mathbb C}/{\mathbb Z}}(M)$ is the subgroup of $\widehat K_{\rm{ch}}(M)$
generated by the triples $(E,\nabla_{\!E},\alpha)$ as above subject to the extra condition that \begin{equation}\label{DankeBunke}d\alpha={\rm{ch}}(\nabla_{\!E})-{\rm{rk}}E
\end{equation}
\subsubsection{Exact sequence and relation with flat bundles:}
The lemma \ref{L:exactseqtransgr} ensures that there is a map from $K^0_{\text{flat}}(M)$ to $\widehat
K_{\rm{ch}}(M)$ which in fact takes its values in $MK^0(M)$ defined by
\begin{equation}\label{platmultiplicative}(E,\nabla_{\!E})\in K^0_{\rm{flat}}(M)\longmapsto(E,\nabla_{\!E},0)
\in\widehat K_{\rm{ch}}(M)
\end{equation}
Then one obtains the following diagram whose lines are exact sequences:
\begin{equation}\begin{CD}K^1_{\text{top}}(M)@>>>K^0_{\text{rel}}(M)@>\partial>>
K^0_{\text{flat}}(M)@>>>K^0_{\text{top}}(M)\\
@VV\Vert V@VV{{\mathcal N}_{\rm{ch}}}V@VVV@VV\Vert V\\
K^1_{\text{top}}(M)@>{\rm{ch}}>>H^{\text{odd}}(M,{\mathbb C})
@>>>MK^0(M)@>{forget}>>K^0_{\text{top}}(M)
\end{CD}\end{equation}
where $\beta\in H^{\rm{odd}}(M,{\mathbb C})$ is sent to
$(E,\nabla_{\!E},\alpha+\beta)-(E,\nabla_{\!E},\alpha)\in
MK^0(M)$ (this element is independent
of the choice of $(E,\nabla_{\!E},\alpha)\in MK^0(M)$ used to represent it).
In this diagram, the part ``$H^{\rm{odd}}(M,{\mathbb C})\longrightarrow MK^0(M)$'' can be replaced
by ``$\Omega^{\rm{odd}}(M,{\mathbb C})\big/d\Omega^{\rm{even}}(M,{\mathbb C})\longrightarrow\widehat K_{\rm{ch}}(M)$'' (defined in the same way) without
losing the commutativity nor the exactness of the second line.

\subsubsection{Borel-type class:}
Put now on any element of $\widehat K_{\rm{ch}}(M)$
\begin{equation}{\mathfrak{B}}(E,\nabla_{\!E},\alpha)=\widetilde{\rm{ch}}
(\nabla_{\!E}^*,\nabla_{\!E})-\alpha+\overline\alpha\end{equation}
where $\nabla_{\!E}^*$ is the adjoint transpose of $\nabla_{\!E}$
for any hermitian metric on $E$. It follows from \eqref{boborelrel} that
\begin{equation}d{\mathfrak{B}}(E,\nabla_{\!E},\alpha)=2i\, {\mathfrak{Im}}\big(\widehat
{\rm{ch}}(E,\nabla_{\!E},\alpha)\big)\end{equation}
Because of relations
\eqref{metricandconnection} and \eqref{Kahatrelation}, this provides a morphism from $\widehat K_{\rm{ch}}(M)$
to the subgroup of purely imaginary forms in $\Omega^{\text{odd}}(M)/d\Omega^{\text{even}}(M)$. It sends $MK^0(M)$ in to $i H^{\rm{odd}}(M,{\mathbb R})$, and from the remark \ref{Lottremark}, one sees that the imaginary part of Cheeger-Chern-Simons classes \cite{CheegerSimons} studied by Bismut and Lott in \cite{BismutLott}
factor throw this morphism and the one defined in
\eqref{platmultiplicative}.

We will denote by $\widehat{{\mathbb R}K}_{\rm{ch}}(M)$ the kernel of ${\mathfrak B}$. It is the invariant subgroup of $\widehat K_{\rm{ch}}(M)$ under the conjugation involution
$(E,\nabla_{\!E},\alpha)\longmapsto(E,\nabla_{\!E}^*,\overline\alpha)$
where $\nabla_{\!E}^*$ is the adjoint transpose of $\nabla_{\!E}$
with respect with any hermitian metric on $E$.
It is exactly the subgroup of $\widehat K_{\rm{ch}}(M)$ generated
by triples $(E,\nabla,\alpha)$ such that $\nabla$ is unitary with respect to some hermitian metric on $E$ and $\alpha$ is real.

It contains the group $K^{-1}_{{\mathbb R}/{\mathbb Z}}(M)$
of $K$-theory with coefficients in ${\mathbb R}/{\mathbb Z}$
considered by Lott in
\cite{Lott} Definition 7. In fact $K^{-1}_{{\mathbb R}/{\mathbb Z}}(M)$ is the intersection of the kernels of $\widehat{\rm{ch}}$ and ${\mathfrak{B}}$.

It follows from lemma \ref{4k+1ou3} that if $E$ is the complexifcation of a real bundle $E_{\!{\mathbb R}}$ on $M$ with connection $\nabla$ coming from a
connection on $E_{\!{\mathbb R}}$, then ${\mathfrak B}(E,\nabla,0)$ vanishes in
degrees $4k+3$ for any integer $k$.
\section{Direct image for relative $K$-theory:}\label{directim}
\subsection{The case of topological $K$-theory:}\label{directimtopologic}
Let $M$ and $B$ be smooth compact real manifolds possibly with
boundary and $\pi\colon M\to B$ a smooth proper submersion.
The fibres of $\pi$ are supposed to be compact without boundary, orientable, and
modelled on the closed manifold $Z$. Moreover the vertical tangent bundle $TZ$, which is the subbundle
of $TM$ consisting of vectors tangent to the fibres of $\pi$, will be
supposed to be globally orientable along $M$. For $y\in B$, $\pi^{-1}(y)$
will be denoted $Z_y$. Let $\xi$ be a ${\mathcal C}^\infty$
complex vector bundle on $M$.
\subsubsection{Construction of the families analytic index map:}\label{s:constr}
For any $y\in B$, consider the infinite dimensional ${\mathbb Z}_2$-graded vector space
\[{\mathcal E}_y^{\pm}
={\mathcal C}^\infty\left(Z_y,\wedge^{\genfrac{}{}{0pt}{}{\text{even}}
{\text{odd}}}T^*Z\otimes\xi\right)=\Omega^{\genfrac{}{}{0pt}{}{\text{even}}
{\text{odd}}}(Z_y,\xi)\]
($T^*Z$ is the dual vector bundle of $TZ$), these spaces are the
fibres of infinite rank vector bundles
${\mathcal E}^+$ and ${\mathcal E}^-$ on $B$ (see \cite{BismutLott} (3.1) to (3.6)) such that
\begin{equation}\label{E:infinitedimvectbundle}
{\mathcal C}^\infty(B,{\mathcal E}^{\pm})
\cong{\mathcal C}^\infty\left(M,\wedge^{\genfrac{}{}{0pt}{}{\text{even}}
{\text{odd}}}T^*Z\otimes\xi\right)
\end{equation}
Choose some connection $\nabla_{\!\xi}$ on $\xi$, then a
vertical exterior differential operator can
be associated to $\nabla_{\!\xi}$:
\[
d^{\nabla_{\!\xi}}\colon\Omega^\bullet(Z_y,\xi)
\longrightarrow \Omega^{\bullet+1}(Z_y,\xi)
\]
$d^{\nabla_{\!\xi}}$ will be considered as an odd endomorphism of the
${\mathbb Z}_2$-graded vector bundle ${\mathcal E}=
{\mathcal E}^+\oplus {\mathcal E}^-$ on $B$ (note that only the restriction of $\nabla_{\!\xi}$ to the fibres of $\pi$ is used to construct $d^{\nabla_{\!\xi}}$).

Choose some smooth hermitian metric $h^\xi$ on $\xi$ and
some smooth euclidean metric $g^Z$ on $TZ$. Associated to $g^Z$ is
a volume form $d{\text{Vol}}_Z$ on the fibres of $\pi$. One
also obtains an associated inner product $(\ /\ )_Z$ on
$\wedge^\bullet T^*Z\otimes\xi$ and the $L^2$ scalar product
on ${\mathcal E}$ defined by
\begin{equation}\label{L2glob}
\langle\alpha,\beta\rangle_{L^2}=\int_{Z_y}
(\alpha/\beta)_Zd{\text{Vol}}_Z\end{equation}
(here $\alpha$ and $\beta\in{\mathcal E}^+_y\oplus{\mathcal E}^-_y$).
This yields a smooth hermitian metric on ${\mathcal E}^+\oplus{\mathcal
E}^-$.
Let $(d^{\nabla_{\!\xi}})^*$ be the formal adjoint of $d^{\nabla_{\!\xi}}$ for this metric.

Let $\eta^+$ and $\eta^-$ be complex vector bundles on $B$ with hermitian metrics $h^+$ and $h^-$.
For
any bundle map $\psi\colon({\mathcal E}^+\oplus\eta^+)\to
({\mathcal E}^-\oplus\eta^-)$ of everywhere finite rank, call
$\psi^*$ the adjoint of $\psi$
with respect to $h^\pm$ and $\langle\ ,\ \rangle_{L^2}$, and set
\begin{align*}
{\mathcal D}^{{\nabla_{\!\xi}}+}_{\psi}&=\left(d^{\nabla_{\!\xi}}+(d^{\nabla_{\!\xi}})^*\right)+\psi
\colon{\mathcal E}^+\oplus\eta^+\longrightarrow
{\mathcal E}^-\oplus\eta^-\\
{\mathcal D}^{{\nabla_{\!\xi}}-}_{\psi}&=\left(d^{\nabla_{\!\xi}}+(d^{\nabla_{\!\xi}})^*\right)
+\psi^*
\colon{\mathcal E}^-\oplus\eta^-\longrightarrow{\mathcal E}^+
\oplus\eta^+
\end{align*}

These are elliptic operators on $Z_y$ so that their kernels are
finite dimensional. A triple $(\eta^+,\eta^-,\psi)$ as above such that
${\text{dim}}\, {\text{Ker}}{\mathcal D}^{{\nabla_{\!\xi}}\pm}_\psi$ are constant
(independent of $y\in B$) will be called
``suitable'' in the sequel.
From \cite{AtiyahSinger} proposition 2.2 (see also \cite{BerlineGetzlerVergne} lemma 9.30 or \cite{SpinGeometry} lemma 8.4 of chapter III),
$\eta^-$ can be chosen null, $\eta^+$ can be chosen as the trivial rank
$N$ vector bundle ${\mathbb C}^N$ (for some great enough $N$) with
some $\psi\colon{\mathbb C}^N\to{\mathcal E}^-$ such that ${\mathcal D}^{{\nabla_{\!\xi}}+}_\psi$ is a surjective
bundle map. The rigidity of the index ensures that
$({\mathbb C}^N,\{0\},\psi)$
is suitable.

For suitable data $(\eta^+,\eta^-,\psi)$, the kernels of ${\mathcal D}^{{\nabla_{\!\xi}}\pm}_\psi$
are vector bundles ${\mathcal H}^\pm$ on $B$.
The class $[{\mathcal H}^+\oplus\eta^-]-[{\mathcal H}^-\oplus\eta^+]\in K_0^{\text{top}}(B)$
co\"\i ncides with $\pi^{\text{Eu}}_*[\xi]$ where the direct image $\pi_*^{\text{Eu}}\colon
K_0^{\text{top}}(M)\to K_0^{\text{top}}(B)$ is
associated to the fibral Euler operator (see \cite{AtiyahSinger} definition 2.3: let $d^Z$ be the usual exterior differential along the fibres acting on usual (vertical) differential form
and $d^{Z*}$ its adjoint with respect to $\langle\ ,\ \rangle$, the fibral Euler operator is $d^Z+d^{Z*}$ acting on vertical differential forms ${\mathbb Z}_2$-graded by the parity of their degree). This
class is independent
on the choice of $\nabla_{\!\xi}$, $h^\xi$ and $g^Z$, it is also independent on the choice of suitable $\eta^\pm$, $h^\pm$ and $\psi$ as will be proved
in the next paragraph. In the particular case above ($({\mathbb C}^N, \{0\},\psi)$
such that ${\mathcal D}^{{\nabla_{\!\xi}}+}_\psi$ is surjective) a proof can be found in
\cite{AtiyahSinger} proposition 2.2 or \cite{BerlineGetzlerVergne} proposition 9.31 (or again \cite{SpinGeometry} lemma 8.4 of chapter III).

Another particular exemple is provided by Mischchenko and
Fomenko's work \cite{MischchenkoFomenko}: it is proved in \cite{MischchenkoFomenko}
that there exists finite rank subbundles ${\mathcal F}^+$ and ${\mathcal F}^-$
of ${\mathcal E}^+$ and ${\mathcal E}^-$ respectively such that $d^{\nabla_{\!\xi}}+(d^{\nabla_{\!\xi}})^*$
is diagonal with respect to the decompositions ${\mathcal E}^\pm={\mathcal F}^\pm
\oplus ({\mathcal F}^\pm)^\perp$ and is invertible on $({\mathcal F}^\pm)^\perp$.
The topological $K$-theoretic direct image is then $[{\mathcal F}^+]-[{\mathcal F}^-]=\pi_*^{\rm{Eu}}[\xi]\, \in\, K^0_{\rm{top}}(B)$. Let $P^{\mathcal F}$
be the orthogonal projection of ${\mathcal E}$ onto ${\mathcal F}$, this
construction corresponds to the choice $\eta^+=\eta^-=\{0\}$ and
$\psi=-P^{\mathcal F}\big(d^{\nabla_{\!\xi}}+(d^{\nabla_{\!\xi}})^*\big)$, and of course ${\mathcal H}^\pm={\mathcal F}^\pm$. 

This construction is functorial by pull-backs on fibered products, which means for maps as:
\begin{equation}\label{pullbackfunctorial}
\begin{CD}\widetilde B\times_{\!B}M@>>>M\\@VVV@VVV\\\widetilde B@>>>B\end{CD}
\end{equation}
(the model of the fibre may of course not change).
\subsubsection{Canonical link between different constructions:}\label{clink}
\underline{Link with representatives obtained from a surjective ${\mathcal D}^+$:}

Choose vector bundles $\eta^+$ and $\eta^-$ on $B$ and $\psi\colon({\mathcal E}^+\oplus\eta^+)\longrightarrow({\mathcal E}^-\oplus\eta^-)$ of finite rank.
With no assumption on the kernels of ${\mathcal D}^{\nabla_{\!\xi}\pm}_\psi$,
it is possible to find some vector bundle $\lambda$ on $B$ and some $\varphi\colon\lambda\to{\mathcal E}^-\oplus\eta^-$ such that ${\mathcal D}^{\nabla_{\!\xi}+}_{\psi+\varphi}$ is surjective
as can be proved in exactly the same way as in \cite{AtiyahSinger} proposition 2.2, or \cite{BerlineGetzlerVergne} lemma 9.30 or \cite{SpinGeometry} lemma 8.4 of chapter III.

Suppose now that $(\eta^+,\eta^-,\psi)$ are suitable and give rise to ${\mathcal H}^\pm$, choose $\lambda$ and $\varphi$ as above. One obtains the following exact sequence of complex vector bundles on $B$:
\begin{equation}\label{E:exact}\begin{aligned}
0\longrightarrow{\mathcal H}^+&\longrightarrow{\text{Ker}}{\mathcal D}^{\nabla_{\!\xi}+}
_{\psi+\varphi}\!\longrightarrow\lambda
\mathop{-\!\!\!-\!\!\!-\!\!\!-\!\!\!\longrightarrow}\limits^{
P^{{\mathcal H}^{\!-}}\!\!\circ\varphi}{\mathcal H}^-\longrightarrow 0\\
&\qquad\ \, (\sigma,v,w)\, \longmapsto w
\end{aligned}\end{equation}
where $P^{{\mathcal H}^-}$ is the projector from ${\mathcal E}^-\oplus\eta^-$ onto
${\mathcal H}^-$ with kernel ${\text{Im}}{\mathcal D}^{\nabla_{\!\xi}+}_{\psi}$.
(In the diagram, ${\text{Ker}}{\mathcal D}^{\nabla_{\!\xi}+}_{\psi+\varphi}$ is seen
as a subbundle of ${\mathcal E}^+\oplus\eta^+\oplus\lambda$
on which its elements are decomposed).

The restriction of ${\mathcal D}^{\nabla_{\!\xi}+}_{\psi+\varphi}$ to
${\mathcal E}^+\oplus\eta^+$ co\"\i ncides with ${\mathcal D}^{\nabla_{\!\xi}+}_\psi$,
so that the obvious inclusion ${\mathcal E}^+\oplus\eta^+\subset {\mathcal E}^+\oplus\eta^+\oplus\lambda$ induces the map ${\mathcal H}^+\subset
{\text{Ker}}{\mathcal D}^{\nabla_{\!\xi}+}
_{\psi+\varphi}$ and the exactness at ${\text{Ker}}{\mathcal D}^{\nabla_{\!\xi}+}
_{\psi+\varphi}$ is obvious. The restriction of ${\mathcal D}^{\nabla_{\!\xi}+}_{\psi+\varphi}$ to
$\lambda$ coincides with $\varphi$,
the surjectivity of ${\mathcal D}^{\nabla_{\!\xi}+}_{\psi+\varphi}$ thus
ensures the surjectivity of
$P^{{\mathcal H}^-}\!\!\circ\varphi$ from $\lambda$ onto
${\mathcal H}^-$. If
$w$ has vanishing image by $P^{{\mathcal H}^-}\!\!\circ\varphi$
then there exists $(\sigma,v)\in{\mathcal E}^+\oplus\eta^+$ such that
${\mathcal D}^{\nabla_{\!\xi}+}_{\psi}(\sigma,v)=\varphi(w)$. Thus
$(-\sigma,-v,w)\in{\text{Ker}}
{\mathcal D}^{\nabla_{\!\xi}+}_{\psi+\varphi}$ and the exactness at
$\lambda$ is proved.

A link between ${\mathcal H}^+-{\mathcal H}^-$ and ${\text{Ker}}{\mathcal D}^{\nabla_{\!\xi}+}
_{\psi+\varphi}-\lambda$ is thus obtained, from which a link between
$({\mathcal H}^+\oplus\eta^-)-({\mathcal H}^-\oplus\eta^+)$ and
$({\text{Ker}}{\mathcal D}^{\nabla_{\!\xi}+}
_{\psi+\varphi}\oplus\eta^-)-(\lambda\oplus\eta^+)$ is trivially deduced.

\noindent\underline{Mutual compatibility of these links:}

Suppose that $\lambda'$ and $\varphi'$ are data satisfying the same surjectivity
hypothesis as $\lambda$ and $\varphi$ with respect to $\eta^\pm$ and $\psi$.
Then $\lambda\oplus\lambda'$ and $\varphi\oplus\varphi'$ also do.
The same construction can be performed from the data
$(\eta^+\oplus\lambda, \eta^-,\psi+\varphi)$ using $\lambda'$ and $\varphi'$
or from the data $(\eta^+\oplus\lambda', \eta^-,\psi+\varphi')$ using
$\lambda$ and $\varphi$. One obtains links between two of the couples
$({\mathcal H}^+\oplus\eta^-)-({\mathcal H}^-\oplus\eta^+)$,
$({\text{Ker}}{\mathcal D}^{\nabla_{\!\xi}+}
_{\psi+\varphi}\oplus\eta^-)-(\lambda\oplus\eta^+)$,
$({\text{Ker}}{\mathcal D}^{\nabla_{\!\xi}+}
_{\psi+\varphi'}\oplus\eta^-)-(\lambda'\oplus\eta^+)$ or
$({\text{Ker}}{\mathcal D}^{\nabla_{\!\xi}+}
_{\psi+\varphi+\varphi'}\oplus\eta^-)-(\lambda\oplus\lambda'\oplus\eta^+)$.
These links are all compatible as can be seen by considering
the exact sequence
\eqref{E:exact} associated to the data $(\lambda\oplus\lambda',\varphi+t\varphi')$
with $t\in[0,1]$ or $(\lambda\oplus\lambda',s\varphi+\varphi')$ with $s\in[0,1]$.

\noindent\underline{How to link any representatives obtained by any suitable data:}

Consider now two suitable data $(\eta^+_0,\eta^-_0,\psi_0)$
and $(\eta_1^+,\eta_1^-,\psi_1)$ which produce ${\mathcal H}_0^\pm$ and ${\mathcal H}_1^\pm$. Call $\eta^+=\eta_0^+\oplus\eta_1^+$,
$\eta^-=\eta_0^-\oplus\eta^-_1$, and $\widetilde{\mathcal E}^\pm$, $\widetilde\eta^\pm$ the pullbacks of ${\mathcal E}^\pm$, $\eta^\pm$
to $\widetilde B=B\times[0,1]$ (by the projection on the first factor). Find some $\widetilde\psi\colon(\widetilde{\mathcal E}^+\oplus\widetilde\eta^+)\longrightarrow
(\widetilde{\mathcal E}^-\oplus\widetilde\eta^-)$
of finite rank which has the expected restrictions $\widetilde\psi\vert_{B\times\{0\}}=\psi_0$ (extended by $0$ on $\eta_1^+$) and $\widetilde\psi\vert_{B\times\{1\}}=\psi_1$ (extended by $0$ on $\eta_0^+$).
The preceding construction (of $\lambda$ and $\varphi$) applies here (on $\widetilde B=B\times[0,1]$) and one
obtains some vector bundle $\widetilde\lambda$ and some bundle map
$\widetilde\varphi\colon\widetilde\lambda\longrightarrow(\widetilde{\mathcal E}^-\oplus\widetilde\eta^-)$ such that
${\mathcal D}^{\nabla_{\!\xi}+}_{\widetilde\psi+\widetilde\varphi}$ is surjective
at any point of $\widetilde B$. Call $\lambda_0=\widetilde\lambda\vert_{B\times\{0\}}$ and similarly for $\lambda_1$,
$\varphi_0=\widetilde\varphi\vert_{B\times\{0\}}$ and $\varphi_1$, denote by $\widetilde{\mathcal G}$ the vector bundle on $\widetilde B$ consisting of the kernel of ${\mathcal D}^{\nabla_{\!\xi}+}_{\widetilde\psi+\widetilde\varphi}$ and ${\mathcal G}_0$ and ${\mathcal G}_1$ its restrictions to $B\times\{0\}$ and $B\times\{1\}$.

The parallel transport along $[0,1]$ produces isotopy classes of isomorphisms
$\lambda_0\overset\sim\longrightarrow\lambda_1$ and ${\mathcal G}_0\overset\sim\longrightarrow{\mathcal G}_1$, from which a link between
$({\mathcal G}_0\oplus\eta^-)-(\lambda_0\oplus\eta^+)$ and $({\mathcal G}_1
\oplus\eta^-)-(\lambda_1\oplus\eta^+)$ is trivially deduced. On $B\times\{0\}$,
$\psi_0$ was extended by $0$ on $\eta_1^+$ to obtain $\psi_0\colon({\mathcal E}^+\oplus\eta_0^+\oplus\eta_1^+)\longrightarrow({\mathcal E}^-\oplus\eta_0^-\oplus\eta_1^-)$, and the kernels to consider become ${\mathcal H}_0^+\oplus\eta_1^+$ and ${\mathcal H}_0^-\oplus\eta^-_1$. The above construction of $\widetilde\lambda$ and $\widetilde\varphi$ restricts then to $B\times\{0\}$ so that a link between
$({\mathcal H}_0^+\oplus\eta_1^+\oplus\eta^-)-({\mathcal H}^-_0\oplus\eta_1^-\oplus\eta^+)$ and
$({\mathcal G}_0\oplus\eta^-)-(\lambda_0\oplus\eta^+)$ is obtained, which is immediately simplified to a link between $({\mathcal H}_0^+\oplus\eta_0^-)-({\mathcal H}^-_0\oplus\eta_0^+)$ and
$({\mathcal G}_0\oplus\eta^-)-(\lambda_0\oplus\eta^+)$. In the same way,
a link is obtained from the restriction to $B\times\{1\}$ between $({\mathcal H}_1^+\oplus\eta_1^-)-({\mathcal H}^-_1\oplus\eta_1^+)$ and
$({\mathcal G}_1\oplus\eta^-)-(\lambda_1\oplus\eta^+)$.
The composition of these three links provide the desired link between
$({\mathcal H}_0^+\oplus\eta_0^-)-({\mathcal H}^-_0\oplus\eta_0^+)$ and
$({\mathcal H}_1^+\oplus\eta_1^-)-({\mathcal H}^-_1\oplus\eta_1^+)$.

\noindent\underline{Independence on the choices and compatibility with composition of links:}

The independence on the choice of $\widetilde\lambda$ and $\widetilde\varphi$
can be shown as was made precedingly for the compatibility of links produced by different choices of $\lambda$ and $\varphi$. The independence on the choice of
$\widetilde\psi$ can be proved by deforming it to any other choice, and make
the above construction on $B\times[0,1]\times[0,1]$.
For three suitable datas $(\eta_i^\pm,\psi_i)$ with $i=1,2,3$,
the same construction can be performed on the product $B\times\Sigma^2$
of $B$ with the two dimensional simplex $\Sigma^2$, and one looks for some
$\widetilde\psi\colon({\mathcal E}^+\oplus\eta_1^+\oplus\eta_2^+\oplus\eta_3^+)
\longrightarrow({\mathcal E}^-\oplus\eta_1^-\oplus\eta_2^-\oplus\eta_3^-)$
with good restrictions on the vertices of $\Sigma^2$ and for some
$\widetilde\lambda$ on $B\times\Sigma^2$ and some $\widetilde\varphi\colon\widetilde\lambda\longrightarrow
({\mathcal E}^-\oplus\eta_1^-\oplus\eta_2^-\oplus\eta_3^-)$ such that
${\mathcal D}^{\nabla_{\!\xi}+}_{\widetilde\psi+\widetilde\varphi}$ be surjective.
The possibility to follow along $\Sigma^2$ the deformation of one edge
to the composition of the two others proves that the constructed links are
mutually compatible.

\noindent\underline{Link between constructions obtained from linked initial data:}

If $\xi^+$, $\xi^-$, $\zeta^+$ and $\zeta^-$ are smooth complex vector bundles on $M$
such that $[\xi^+]-[\xi^-]=[\zeta^+]-[\zeta^-]\in K^0_{\text{top}}(M)$,
choose connections $\nabla_{\!\xi^+}$, $\nabla_{\!\xi^-}$, $\nabla_{\!\zeta^+}$
and $\nabla_{\!\zeta^-}$ and consider some vector bundle isomorphism $\ell\colon
\xi^+\oplus\zeta^-\oplus L\overset\sim\longrightarrow\xi^-\oplus\zeta^+\oplus
L$ as in \eqref{linque}. Choose some connection $\nabla_{\!L}$ on $L$.
Pull-back the vector bundle $\xi^+\oplus\zeta^-\oplus L$
on $M\times[0,1]$ and endow this pull-back with some connection $\widetilde\nabla$
such that its restrictions to $M\times\{0\}$ and $M\times\{1\}$ equal
$\nabla_{\!\xi^+}\oplus\nabla_{\!\zeta^-}\oplus\nabla_{\!L}$
and $\ell^*(\nabla_{\!\xi^-}\oplus\nabla_{\!\zeta^+}\oplus\nabla_{\!L})$
respectively. Then the construction of analytic families index for the
submersion $\pi\times{\rm{Id}}_{[0,1]}\colon M\times[0,1]\longrightarrow B\times[0,1]$ restricts on $M\times\{0\}$ and $M\times\{1\}$
to regular constructions of analytic families index representatives
of $\pi_*^{\rm{Eu}}[\xi^+]+\pi_*^{\rm{Eu}}[\zeta^-]+\pi_*^{\rm{Eu}}
[L]$ and $\pi_*^{\rm{Eu}}[\xi^-]+\pi_*^{\rm{Eu}}[\zeta^+]+\pi_*[L]$ respectively.
The parallel transport along $[0,1]$ thus provides a link between direct images
of $\xi^+-\xi^-$ and $\zeta^+-\zeta^-$ (in a little more complicated form than in \eqref{linque}).

The independency of the class of this link on the way of constructing the families analytic index for the submersion $\pi\times{\rm{Id}}_{[0,1]}$ is immediate, the independency on the choice of $L$ and $\ell$ (in some same equivalence class of links)
can be proved by deforming the link as in \eqref{linkdelink} and constructing a
direct image for a submersion of the form $\pi\times{\rm{Id}}_{[0,1]\times[0,1]}$.

\noindent\underline{The case of subspaces of ${\mathcal E}$ respected by $d^{\nabla_{\!\xi}}+(d^{\nabla_{\!\xi}})^*$:}

Suppose that ${\mathcal H}^\pm\subset{\mathcal F}^\pm\subset{\mathcal E}^\pm$
are finite rank subbundles of ${\mathcal E}^\pm$ respected by $d^{\nabla_{\!\xi}}
+(d^{\nabla_{\!\xi}})^*$ and such that $d^{\nabla_{\!\xi}}+(d^{\nabla_{\!\xi}})^*$ is invertible on ${\mathcal H}^\perp$ and also on ${\mathcal F}^\perp$. Call ${\mathcal K}^\pm$ the
orthocomplement of ${\mathcal H}^\pm$ in ${\mathcal F}^\pm$.
$d^{\nabla_{\!\xi}}+(d^{\nabla_{\!\xi}})^*$ respects ${\mathcal K}^\pm$ and is
invertible on it.
Let $P^{{\mathcal K}^\pm}$, $P^{{\mathcal H}^\pm}$ and $P^{{\mathcal F}^\pm}$ be the orthogonal projections onto
${\mathcal K}^\pm$, ${\mathcal H}^\pm$ and ${\mathcal F}^\pm$ respectively.
Make the construction of the third alinea above with $\widetilde\psi=-(P^{{\mathcal H}^-}+tP^{{\mathcal K}^-})\big(d^{\nabla_{\!\xi}}+(d^{\nabla_{\!\xi}})^*\big) P^{{\mathcal F}^+}$, (so that ${\mathcal H}^\pm$
correspond to $t=0$ and ${\mathcal F}^\pm$ to $t=1$), and take
$\widetilde\lambda={\mathcal F}^-$ and $\widetilde\varphi$ the natural inclusion
of ${\mathcal F}^-$ into ${\mathcal E}^-$. Call as above $\widetilde{\mathcal G}$ the resulting kernel of ${\mathcal D}^{\nabla_{\!\xi}+}_{\widetilde\psi+\widetilde\varphi}$ on $B\times[0,1]$.

The canonical link on $B\times\{1\}$ is trivial; for any $t$, ${\mathcal G}_t$ is the image of
${\mathcal F}^+$ in ${\mathcal F}^+\oplus{\mathcal F}^-$
by ${\rm{Id}}_{{\mathcal F}^+}\oplus(t-1)P^{{\mathcal K}^-}\big(d^{\nabla_{\!\xi}}+(d^{\nabla_{\!\xi}})^*\big)$, so that
the parallel transport along $[0,1]$ induces an isomorphism between ${\mathcal G}_0$ and ${\mathcal F}^+$. Throw this isomorphism,
the exact sequence \eqref{E:exact} describing the canonical link on $B\times\{0\}$ reads as \eqref{oscour}:
\begin{equation}\label{linkvaleurspropres}
0\longrightarrow{\mathcal H}^+\longrightarrow{\mathcal F}^+
\overset{-P^{{\mathcal K}^-}\big(d^{\nabla_{\!\xi}}+(d^{\nabla_{\!\xi}})^*\big)}{-\!\!\!-\!\!\!-\!\!\!-\!\!\!-\!\!\!-\!\!\!-\!\!\!-\!\!\!-\!\!\!-\!\!\!-\!\!\!-\!\!\!\longrightarrow}{\mathcal F}^-\longrightarrow{\mathcal H}^-
\longrightarrow0\end{equation}
The canonical link constructed precedingly between
${\mathcal H}^+-{\mathcal H}^-$ and ${\mathcal F}^+-{\mathcal F}^-$
is nothing but the link obtained as in \eqref{precedente}
from this exact sequence.
\subsection{Fibral Hodge symmetry:}
\subsubsection{The fibral Hodge operator:}
For any vertical tangent vector ${\tt w}\in TZ$, consider its dual
one-form ${\tt w}^\flat$ (throw the metric $g^Z$), and its Clifford action \begin{equation}\label{Clifford}
c({\tt w})=({\tt w}^\flat\wedge)-\iota_{\tt w}
\end{equation}
on $\wedge^\bullet T^*Z\otimes\xi$ ($\iota_{\tt w}$
denotes the interior product by ${\tt w}$); $c({\tt w})$ is skewadjoint
with respect to $(\ /\ )_Z$ and verifies $c({\tt w})^2=-g^Z({\tt w},{\tt w})$,
it is an isometry with respect to $(\ /\ )_Z$ if $g^Z({\tt w},{\tt w})=1$.

Consider the vertical Hodge operator $*_Z=c({\tt e}_1)c({\tt e}_2)\ldots c({\tt e}_{{\rm{dim}}Z})$ for any orthonormal direct base
${\tt e}_1,{\tt e}_2,\ldots,{\tt e}_{{\rm{dim}}Z}$ of $TZ$. It is an isometry
of ${\mathcal E}$ (endowed with $\langle\ ,\ \rangle_{L^2}$), and it has the same parity as ${\rm{dim}}Z$ (with respect to the ${\mathbb Z}_2$ grading of ${\mathcal E}$). Its inverse $*_Z^{-1}=(-1)^{\frac12{\rm{dim}}Z({\rm{dim}}Z+1)}*_Z$ is also its adjoint
with respect to both $(\ /\ )_Z$ and $\langle\ ,\ \rangle_{L^2}$. Define the metrized exterior product of $\xi$-valued
vertical differential forms by the following formula on decomposed tensors:
\[(\alpha\widehat\otimes a)\mathop{\wedge}\limits_{h^\xi}(\beta\widehat\otimes b)
=(\alpha\wedge\overline\beta) h^\xi(a,b)\]
(a sign $(-1)^{{\rm{deg}}a{\rm{deg}}\beta}$ should be put on the right side
if $\xi$ would be ${\mathbb Z}_2$-graded, but this case will not be considered in the sequel, note also that this operation is independent of the riemannian vertical metric $g^Z$).
Then for any $\gamma\in{\mathcal E}$ whose differential form degree is $\leq{\rm{deg}}\alpha$:
\begin{equation}\label{starwedge}(\alpha\widehat\otimes a)\mathop{\wedge}\limits_{h^\xi}(*_Z\gamma)=(
-1)^{\frac12{\rm{deg}}\alpha({\rm{deg}}\alpha-1)
+{\rm{dim}}Z{\rm{deg}}\alpha}\big((\alpha\widehat\otimes a)\big/\gamma\big)_Zd{\rm{Vol}}_Z\end{equation}

For any ${\tt w}\in TZ$, $c({\tt w})$ commutes with $*_Z$ if ${\rm{dim}}Z$ is odd and
it anticommutes with $*_Z$ if ${\rm{dim}}Z$ is even. It follows from the two preceding formulae that
if $\nabla_{\!\xi}^*$ is associated to $\nabla_{\!\xi}$ and $h^\xi$
as in \eqref{adjointconnection}, then for any $\gamma$ and $\gamma'$ in
${\mathcal E}$:
\begin{equation}\label{flipflop}
\begin{aligned}
d^Z(\gamma\underset{h^\xi}\wedge\gamma')&=(d^{\nabla_{\!\xi}}\gamma)\underset{h^\xi}\wedge\gamma'
+(-1)^{{\rm{deg}}\gamma}\gamma\underset{h^\xi}\wedge(d^{\nabla_{\!\xi}^*}\gamma')\\
{\text{so that}}\qquad (d^{\nabla_{\!\xi}})^*&=(-1)^{1+\frac12{\rm{dim}}Z({\rm{dim}}Z-1)}*_Zd^{\nabla_{\!\xi}^*}*_Z
\end{aligned}\end{equation}
from which one deduces that
\begin{equation}\label{agentdechange}
d^{\nabla_{\!\xi}}+(d^{\nabla_{\!\xi}})^*=-(-1)^{{\rm{dim}}Z}*_Z^{-1}
\big(d^{\nabla_{\!\xi}^*}+(d^{\nabla_{\!\xi}^*})^*\big)\, *_Z
\end{equation}
This formula
can also be checked
from \cite{BismutLott} formulae (3.36), (1.30), (1.31) and the last sentence at the end of the first alinea of \S III(d).
\subsubsection{Fibral Hodge symmetry and families index: the even dimensional fibre case}\label{evenHodge}
Suppose that ${\rm{dim}}Z$ is even. Endow $\xi$ with some hermitian metric $h^\xi$ and connection $\nabla_{\!\xi}$, choose suitable data
$\eta^+$, $\eta^-$ and $\psi$ to construct a representative $[{\mathcal H}^+\oplus\eta^-]-[{\mathcal H}^-\oplus\eta^+]$ of $\pi^{\rm{Eu}}_*[\xi]\in K^0_{\rm{top}}(B)$. It is then deduced from \eqref{agentdechange} that the data $\big(\eta^+,\eta^-,-(*_Z\oplus{\rm{Id}}_{\eta^-})\circ\psi\circ
(*_Z^{-1}\oplus{\rm{Id}}_{\eta^+})\big)$ are suitable with respect to $\xi$ endowed with $h^\xi$ and $\nabla^*_{\!\xi}$, and that they
produce the following representative $[(*_Z\oplus{\rm{Id}}_{\eta^+}){\mathcal H}^+\oplus\eta^-]
-[(*_Z\oplus{\rm{Id}}_{\eta^-}){\mathcal H}^-\oplus\eta^+]$
of $\pi^{\rm{Eu}}_*[\xi]$.

Thus $*_Z$ gives an
isomorphism (and thus an equivalence class of link which will be denoted $[\ell_*]$) between representatives of $\pi_*^{\rm{Eu}}[\xi]$ obtained
from $\nabla_{\!\xi}$ and $\nabla_{\!\xi}^*$ (in this order, say).
And the generic
construction of paragraph \ref{clink} also produces some canonical equivalence class of link $[\ell_{\rm{can}}]$ between them.
\begin{lemma}\label{starlink} $[\ell_*]=[\ell_{\rm{can}}]$
\end{lemma}
\begin{proof} Put $\nabla_{\!\xi}^u=\frac12(\nabla_{\!\xi}+\nabla_{\!\xi}^*)$
and consider some suitable data $(\lambda,\{0\},\varphi)$ with respect to
$\xi$ endowed with $h^\xi$ and $\nabla_{\!\xi}^u$, such that
$\varphi$ vanishes on ${{\mathcal E}^+}$ and ${\mathcal D}^{\nabla_{\!\xi}^u+}_\varphi$ is surjective. (Such data exist as
was used at the beginning of \S\ref{clink}). As just seen above, $(\lambda,\{0\},\Phi)$ with $\Phi=-*_Z\circ\varphi$ is also suitable, and of course
${\mathcal D}^{\nabla_{\!\xi}^u+}_\Phi$ is also surjective.

The data
$(\lambda\oplus\lambda,\{0\},\varphi+\Phi)$ (where $\varphi$ only
acts on the first copy of $\lambda$ and $\Phi$ only on the second one) are also suitable (as was used in the second alinea of \S\ref{clink}), and ${\mathcal D}_{\varphi+\Phi}^{\nabla_{\!\xi}^u+}$
is surjective. Call ${\mathcal H}^{\frac12}$ its kernel and $\varepsilon=(-1)^{
\frac12{\rm{dim}}Z({\rm{dim}}Z+1)}$. ${\mathcal H}^{\frac12}$ is, as subbundle of ${\mathcal E}^+\oplus\lambda\oplus\lambda$, invariant under $\sigma_{\mathcal H}=*_Z\oplus(^0_1{\, }^\varepsilon_0)$, where
$(^0_1{\, }^\varepsilon_0)$ is to understand as $(v,w)\in\lambda\oplus\lambda\longmapsto(\varepsilon w,v)\in\lambda\oplus\lambda$.
The canonical equivalence classes of links $[\ell_\varphi]$ between $({\rm{Ker}}{\mathcal D}^{\nabla_{\!
\xi}^u}_\varphi)-\lambda$ and ${\mathcal H}^{\frac12}-(\lambda\oplus\lambda)$, and $[\ell_\Phi]$ between $({\rm{Ker}}{\mathcal D}^{\nabla_{\!
\xi}^u}_\Phi)-\lambda$ and ${\mathcal H}^{\frac12}-(\lambda\oplus\lambda)$ are conjugated in the sense that
\begin{equation}\label{relconnect}[\ell_\varphi]=\big[\sigma_{\mathcal H}\oplus(^0_1\, ^\varepsilon_0)\big]\circ[\ell_\Phi]\circ\big[(*_Z\oplus{\rm{Id}}_{\lambda})\oplus{\rm{Id}}_{\lambda}\big]
\end{equation}

Choose some path $(\nabla_{\!t})_{t\in[0,1]}$ of connections on $\xi$
such that $\nabla_{\!0}=\nabla_{\!\xi}$ and $\nabla_{\!1}=\nabla^u_{\!\xi}$.
On $M\times[0,1]$, pull back $\xi$ with $h^\xi$ and endow its
pullback $\widetilde\xi$ with the connection $\widetilde\nabla=d_{[0,1]}+\nabla_{\!t}$.
Choose suitable data $(\widetilde\nu^+,\widetilde\nu^-,\widetilde\psi)$
with respect to this situation. As previously, the parallel transport
along $[0,1]$ with adapted canonical links on $B\times\{0\}$ and $B\times\{1\}$ provides a model for the canonical link $[\ell_{\nabla_{\!\xi}}]$
between $({\mathcal H}^+\oplus\eta^-)-({\mathcal H}^-\oplus\eta^+)$
and ${\mathcal H}^\frac12-(\lambda\oplus\lambda)$.

Endow now $\widetilde\xi$ with $\widetilde\nabla^*=d_{[0,1]}+
\nabla_{\!t}^*$, and consider the associated suitable data
$\big(\widetilde\nu^+,\widetilde\nu^-,-(*_Z\oplus{\rm{Id}}_{\widetilde\nu^-})\circ\widetilde\psi\circ
(*_Z^{-1}\oplus{\rm{Id}}_{\widetilde\nu^+})\big)$. One obtains a model for the canonical link
$[\ell_{\nabla_{\!\xi}^*}]$ between $\big((*_Z\oplus{\rm{Id}}_{\eta^+}){\mathcal H}^+\oplus\eta^-\big)-\big((*_Z\oplus{\rm{Id}}_{\eta^-}){\mathcal H}^-\oplus\eta^+\big)$
and ${\mathcal H}^\frac12-(\lambda\oplus\lambda)$. In one hand, from the mutual compatibility of canonical links, one has $[\ell_{\rm{can}}]=[\ell_{\nabla_{\!\xi}}\circ\ell_{\nabla_{\!\xi}^*}^{-1}]$.
In the other hand, the fact that these two models are obtained throw $*_Z$ from one another proves that $[\ell_{\nabla_{\!\xi}}]$ and $[\ell_{\nabla^*_{\!\xi}}]$
are related in the same way as $[\ell_\varphi]$ and $[\ell_\Phi]$ are in \eqref{relconnect}:
\begin{equation}\label{varphiPhi}
[\ell_{\nabla_{\!\xi}}]=\big[\sigma_{\mathcal H}\oplus(^0_1\, ^\varepsilon_0)\big]\circ[\ell_{\nabla_{\!\xi}^*}]\circ\big[(*_Z\oplus{\rm{Id}}_{\eta^-})\oplus{\rm{Id}}_{\eta^+}\big]=\big[\sigma_{\mathcal H}\oplus(^0_1\, ^\varepsilon_0)\big]\circ[\ell_{\nabla_{\!\xi}^*}]\circ[\ell_*]\end{equation}

Indeed, the construction of adapted canonical links on $B\times\{0\}$ and $B\times\{1\}$ entering in the construction of $[\ell_{\nabla_{\!\xi}}]$ and
their counterpart in the construction of $[\ell_{\nabla_{\!\xi}^*}]$ are
easily seen to be $*_Z$-conjugated from
the fact that the construction of the exact sequence \eqref{E:exact} before
or after letting $*_Z\oplus{\rm{Id}}_{\eta^+}\oplus{\rm{Id}}_{\eta^-}\oplus{\rm{Id}}_\lambda$ act on its components, gives the same result.

It follows that the obstruction to $[\ell_*]$ being equal to $[\ell_{\rm{can}}]$ would
be the nontriviality of the link between ${\mathcal H}^{\frac12}-(\lambda\oplus\lambda)$ and itself provided by the au\-to\-mor\-phism $\sigma_{\mathcal H}\oplus(^0_1\, ^\varepsilon_0)$, i.e. the nontriviality of the classes in $K^1_{\rm{top}}(B)$ defined by the automorphism $\sigma_{\mathcal H}$ of ${\mathcal H}^\frac12$ and the automorphism $(^0_1\, ^\varepsilon_0)$ of $\lambda\oplus\lambda$.
In any case, the square of these automorphisms is the identity or
minus the identity (on any point of $B$), so that its eigenspaces are
vector bundles on $B$, on which it reduces to a constant multiple of the identity: the obtained element of $K^1_{\rm{top}}(B)$ thus vanishes and the lemma is proved.
\end{proof}
\subsubsection{Fibral Hodge symmetry and families index: the odd dimensional fibre case}\label{oddHodge}
Suppose now that ${\rm{dim}}Z$ is odd, and consider the path $\nabla_{\!t}=(1-t)\nabla_{\!\xi}
+t\nabla_{\!\xi}^*$ of connections on $\xi$, which of course verifies that $\nabla_{\!t}
=\nabla_{\!1-t}^*$ for all $t\in[0,1]$. Pull back $\xi$ on $M\times[0,1]$
and endow the obtained bundle $\widetilde\xi$ on $M\times[0,1]$ with the connection $\widetilde\nabla=d_{[0,1]}+\nabla_{\!t}$ (and $h^\xi$).
Call $\widetilde{\mathcal E}^\pm$ the pullback of ${\mathcal E}^\pm$ on $B\times[0,1]$
and ${\mathcal E}^\pm_t$ its restriction to $B\times\{t\}$.

As proved in \cite{MischchenkoFomenko},
there exists finite rank subbundles $\widetilde{\mathcal F}^+$ and $\widetilde{\mathcal F}^-$
of $\widetilde{\mathcal E}^+$ and $\widetilde{\mathcal E}^-$ respectively such that $d^{\widetilde\nabla}+(d^{\widetilde\nabla})^*$
is diagonal with respect to the decompositions $\widetilde{\mathcal E}^\pm=\widetilde{\mathcal F}^\pm
\oplus (\widetilde{\mathcal F}^\pm)^\perp$ and is invertible on $(\widetilde{\mathcal F}^\pm)^\perp$.

For any $t\in[0,1]$, put ${\mathcal G}_t^+=*_Z{\mathcal F}_{1-t}^-$, and ${\mathcal G}_t^-=*_Z{\mathcal F}_{1-t}^+$.
This yields a subvector bundle $\widetilde{\mathcal G}^\pm$ of $\widetilde{\mathcal E}^\pm$ on $B\times[0,1]$ whose restrictions to
$B\times\{t\}$ are ${\mathcal G}^\pm_t$ for any $t$.
From the relation \eqref{agentdechange}, which is also valid for $d^{\widetilde\nabla}
+(d^{\widetilde\nabla})^*$, it follows that $d^{\widetilde\nabla}+(d^{\widetilde\nabla})^*$
is diagonal with respect to the decompositions $\widetilde{\mathcal E}^\pm=\widetilde{\mathcal G}^\pm
\oplus (\widetilde{\mathcal G}^\pm)^\perp$ and is invertible on $(\widetilde{\mathcal G}^\pm)^\perp$.

Call $\widetilde{\mathcal H}$ the sum of $\widetilde{\mathcal F}$ and $\widetilde{\mathcal G}$ as subbundles of ${\mathcal E}$.
As $\widetilde{\mathcal F}$ and $\widetilde{\mathcal G}$, $\widetilde{\mathcal
H}$ has the property that $d^{\widetilde\nabla}+(d^{\widetilde\nabla})^*$
is invertible on $(\widetilde{\mathcal H}^\pm)^\perp$ and diagonal with respect to the decomposition $\widetilde{\mathcal E}^\pm
=\widetilde{\mathcal H}^\pm\oplus(\widetilde{\mathcal H}^\pm)^\perp$.
Call ${\mathcal H}_t^\pm$ the restrictions of $\widetilde{\mathcal H}^\pm$
to $B\times\{t\}$ for any $t$, then $\widetilde{\mathcal H}$ has the extra property that ${\mathcal H}^+_t=*_Z{\mathcal H}^-_{1-t}$ and ${\mathcal H}^-_t=
*_Z{\mathcal H}^+_{1-t}$ for any $t\in[0,1]$.

Thus $*_Z$ realises isomorphisms $*_Z\colon{\mathcal H}_0^+\overset
\sim\longrightarrow{\mathcal H}^-_1$ and $*_Z\colon{\mathcal H}^-_0
\overset\sim\longrightarrow{\mathcal H}^+_1$. In the other hand, the parallel transport along $[0,1]$
realises isotopy classes of isomorphisms $g^+_\parallel\colon
{\mathcal H}^+_0\overset\sim\longrightarrow{\mathcal H}^+_1$ and
$g^-_\parallel\colon{\mathcal H}^-_0\overset\sim\longrightarrow{\mathcal H}^-_1$. As a first consequence, it follows that ${\mathcal H}^+_0$
and ${\mathcal H}^-_0$ are isomorphic throw $(g^-_\parallel)^{-1}\circ
*_Z$. In fact, as was pointed out in \S\ref{s:constr},
$[{\mathcal H}^+_0]-[{\mathcal H}^-_0]$ is a representative of $\pi^{\rm{Eu}}_*[\xi]$
in $K^0_{\rm{top}}(B)$.
This proves that $\pi^{\rm{Eu}}_*\colon K^0_{\rm{top}}(M)\to K^0_{\rm{top}}(B)$ vanishes in the case of odd dimensional fibres.

The isomorphism between ${\mathcal H}^+_0$ and ${\mathcal H}^-_0$
can be interpreted as a link $[\ell_{{\mathcal H}_0}^{\{0\}}]$ between ${\mathcal H}^+_0$ and
${\mathcal H}^-_0$ with $\{0\}-\{0\}$.
The isomorphism $*_Z\colon{\mathcal H}^+_t\overset\sim\longrightarrow
{\mathcal H}^-_{1-t}$ is continuous with respect to $t$, so that
$(g^-)^{-1}_\parallel\circ{*_Z}=*_Z^{-1}\circ g^+_\parallel$. As was
seen before $*_Z^{-1}=\pm*_Z$, and the change of sign is isotopic to
the identity, so that the link $[\ell_{{\mathcal H}_0}^{\{0\}}]$ can be indifferently be obtained
from $(g^-_\parallel)^{-1}\circ
*_Z$ or from $*_Z\circ g^+_\parallel$.
The corresponding link $[\ell^{\{0\}}_{{\mathcal H}_1}]$
between ${\mathcal H}^+_1-{\mathcal H}^-_1$ and $\{0\}-\{0\}$ obtained
from $*_Z\circ(g^-_\parallel)^{-1}$ or $g^+_\parallel\circ
*_Z$ verifies $[\ell_{{\mathcal H}_1}^{\{0\}}]=*_Z\circ[-\ell_{{\mathcal H}_0}^{\{0\}}]$ (the $-$ sign comes from the fact that $*_Z$ realises a link between ${\mathcal H}^+_1-{\mathcal H}^-_1$
and ${\mathcal H}^-_0-{\mathcal H}^+_0$).

This can be generalised in the following way: let $(\eta^+,\eta^-,\psi)$ be suitable data
for $\xi$ on $M$ endowed with $\nabla_{\!\xi}$ and $h^\xi$.
Call ${\mathcal K}^\pm_0={\rm{Ker}}{\mathcal D}^{\nabla_{\!\xi}\pm}_\psi$
the obtained vector bundles on $B$. The canonical link $[\ell^0_{\rm{can}}]$ between $({\mathcal K}_0^+\oplus\eta^-)-({\mathcal K}_0^-\oplus\eta^+)$ and ${\mathcal H}^+_0-{\mathcal H}^-_0$
and $[\ell_{{\mathcal H}_0}^{\{0\}}]$ compose to provide a link between $({\mathcal K}_0^+\oplus\eta^-)-({\mathcal K}_0^-\oplus\eta^+)$ and $\{0\}-\{0\}$ which will be called $[\ell^{\{0\}}_{{\mathcal K}_0}]$.
\begin{lemma}\label{linkzero}
This link $[\ell^{\{0\}}_{{\mathcal K}_0}]$ is canonical.
\end{lemma}
\begin{proof}
This means that it is independent on the chosen $\widetilde{\mathcal H}$ used
in its construction (in particular, it is compatible with canonical links between different representatives of the direct image).

If $\widetilde{\mathcal H}'$ is another one (respected by $d^{\widetilde\nabla}
+(d^{\widetilde\nabla})^*$, such that $d^{\widetilde\nabla}+(d^{\widetilde\nabla})^*$ is invertible on
its orthocomplement, and such that ${\mathcal H}_t'{}^\pm=*_Z{\mathcal H}'{}^\mp_{\!\!\!1-t}$) then their sum and their intersection also are, so that we may and will suppose that $\widetilde{\mathcal H}'\subset\widetilde{\mathcal H}$. The orthocomplement $\widetilde{\mathcal H}''$ of $\widetilde{\mathcal H}'$
in $\widetilde{\mathcal H}$ is respected by $d^{\widetilde\nabla}+(d^{\widetilde\nabla})^*$,
which is invertible on it, and it has the symmetry property
${\mathcal H}_t''{}^\pm=*_Z{\mathcal H}''{}^\mp_{\!\!\!\!\!1-t}$. Thus the previous isomorphism ${\mathcal H}^+_0\cong{\mathcal H}^-_0$ given by $(g^-_\parallel)^{-1}\circ*_Z$ respects the decomposition ${\mathcal H}={\mathcal H}'\oplus{\mathcal H}''$. Let $P^{{\mathcal H}_t'}$ and
$P^{{\mathcal H}''_t}$ be the orthogonal projections of ${\mathcal H}_t$
onto ${\mathcal H}'_t$ and ${\mathcal H}''_t$ respectively, then
this proves that the link $[\ell_{{\mathcal H}_0}^{\{0\}}]\circ[\ell_{{\mathcal H}'_0}^{\{0\}}]$ between ${\mathcal H}'_0{}^{\!+}-{\mathcal H}'_0{}^{\!-}$
and ${\mathcal H}^+_0-{\mathcal H}^-_0$ (provided by the fact that they are both linked to $\{0\}-\{0\}$), is the same as the link provided as in \eqref{precedente} by the following exact sequence:
\[0\longrightarrow{\mathcal H}'_0{}^{\!+}\overset\subset\longrightarrow
{\mathcal H}^+_0\overset{P^{{\mathcal H}''_0}\big((g_\parallel^-)^{-1}\circ*_Z\big)}
{-\!\!\!-\!\!\!-\!\!\!-\!\!\!-\!\!\!-\!\!\!-\!\!\!-\!\!\!-\!\!\!-\!\!\!\longrightarrow}{\mathcal H}^-_0\overset{P^{{\mathcal H}'_0}}\longrightarrow{\mathcal H}'_0{}^{\!-}
\longrightarrow0\]

In the other hand, the canonical link between ${\mathcal H}'_0{}^{\!+}-{\mathcal H}'_0{}^{\!-}$
and ${\mathcal H}^+_0-{\mathcal H}^-_0$ provided by the fact that they are
both representatives of the direct image $\pi^{\rm{Eu}}_*[\xi]$ is the one
associated as in \eqref{precedente} to the following exact sequence:
\[0\longrightarrow{\mathcal H}'_0{}^{\!+}\overset\subset\longrightarrow
{\mathcal H}^+_0\overset{P^{{\mathcal H}''_0}{\mathcal D}_0}
{-\!\!\!-\!\!\!-\!\!\!-\!\!\!-\!\!\!-\!\!\!-\!\!\!-\!\!\!\longrightarrow}{\mathcal H}^-_0\overset{P^{{\mathcal H}'_0}}\longrightarrow{\mathcal H}'_0{}^{\!-}
\longrightarrow0\]
where ${\mathcal D}_0$ is the restriction of $d^{\widetilde\nabla}+
(d^{\widetilde\nabla})^*$ to $M\times\{0\}$ (see \eqref{linkvaleurspropres}). The lemma thus follows from the triviality of the element of $K^1_{\rm{top}}(B)$ defined by the automorphism ${\mathcal D}_0\circ\big((g_\parallel^-)^{-1}\circ*_Z\big)^{-1}$ of ${\mathcal H}''_0{}^+$ (because of the compatibility of canonical links constructed in \S\ref{clink}).

This automorphism is isotopic by parallel transport to the automorphism ${\mathcal D}_{\frac12}\circ*_Z^{-1}$ of ${\mathcal H}''_{\frac12}{}^+$ (where ${\mathcal D}_{\frac12}$
is the restriction of $d^{\widetilde\nabla}+(d^{\widetilde\nabla})^*$ to
$M\times\{\frac12\}$). And this last one is either selfadjoint or skewadjoint
according to ${\rm{dim}}Z$ (combine \eqref{agentdechange} with the fact that $*_Z^{-1}=(-1)^{\frac12{\rm{dim}}Z({\rm{dim}}Z+1)}*_Z$ is the adjoint of $*_Z$). The conclusion follows the considerations at the end of \S\ref{cathop}.
\end{proof}
Back to ${\mathcal K}^\pm_0$, it follows from
\eqref{agentdechange} that
$\big(\eta^-, \eta^+,(*_Z\oplus{\rm{Id}}_{\eta^+})\circ
\psi^*\circ(*_Z^{-1}\oplus{\rm{Id}}_{\eta^-})\big)$ are suitable for $\xi$ endowed with $h^\xi$ and $\nabla_{\!\xi}^*$, and produce the vector bundles
${\mathcal K}^+_1=(*_Z\oplus{\rm{Id}}_{\eta^-}){\mathcal K}^-_0$ and
${\mathcal K}^-_1=(*_Z\oplus{\rm{Id}}_{\eta^+}){\mathcal K}^+_0$. The link $[\ell_{{\mathcal K}_1}^{\{0\}}]$ between $({\mathcal K}_1^+\oplus\eta^+)-({\mathcal K}_1^-\oplus\eta^-)$ and $\{0\}-\{0\}$ obtained as above for ${\mathcal K}_0^\pm$, throw
${\mathcal H}^+_1-{\mathcal H}_1^-$, verifies the equality $[\ell_{{\mathcal K}_1}^{\{0\}}]
=\big((*_Z\oplus{\rm{Id}}_{\eta^-\oplus\eta^+})\oplus(*_Z\oplus{\rm{Id}}_{\eta^+\oplus\eta^-})\big)\circ[-\ell_{{\mathcal K}_0}^{\{0\}}]$.
The reason for this is that $[\ell^0_{\rm{can}}]$ and $[\ell_{\rm{can}}^1]=[\ell_{{\mathcal H}_1}^{\{0\}}]^{-1}\circ[\ell_{{\mathcal K}_1}^{\{0\}}]$
are $*_Z$-conjugated is the same way as in the even dimensional fibre case
\eqref{varphiPhi} (but with a $-$ sign here).
\subsection{Flat $K$-theory:}

To some flat vector bundle $(E,\nabla_{\!E})$ on $M$, the
de Rham cohomology $H^\bullet(Z,E)$ of the fibres of $\pi$ with coefficients in $E$
provides (graded) vector bundles on $B$, which are endowed with flat connections in a canonical way, see \cite{BismutLott} \S III (f).

Put
$\pi_!^+E=H^{\rm{even}}(Z,E)$ and $\pi_!^-E=H^{\rm{odd}}(Z,E)$,
and call $\nabla_{\!\pi_!^+E}$ and $\nabla_{\!\pi_!^-E}$ their canonical flat connections ($\pi^\pm_!E$ are smooth vector bundles on $B$, but whose definition depends on $\nabla_{\!E}$), then the map \[(E,\nabla_{\!E})\longmapsto(\pi^+_!E,\nabla_{\!\pi_!^+E})-
(\pi^-_!E,\nabla_{\!\pi_!^-E})\]
provides a
morphism $\pi_!\colon K^0_{\text{flat}}
(M)\to K^0_{\text{flat}}(B)$ (see \cite{BismutLott} Proposition 3.14).
\subsubsection{Direct images and short exact sequences:}
For a short exact sequence of flat bundle as in \eqref{suitexacteplate}
\begin{equation}\label{suitexplate}
0\longrightarrow (E',\nabla_{\!E'})\overset i\longrightarrow
(E,\nabla_{\!E})
\overset p\longrightarrow(E'',\nabla_{\!E''})\longrightarrow0
\end{equation}
the long exact sequence in cohomology reads
\begin{equation}\label{diagrexplat}\begin{CD}\pi_!^+E'@>[i]>>\pi^+_!E
@>[p]>>\pi^+_!E''\\@AAA@.@VVV\\\pi^-_!E''
@<[p]<<\pi_!^-E@<[i]<<\pi_!^-E'\end{CD}\end{equation}
this turns out to be an exact sequence of flat vector bundles on $B$ and this proves that $\pi_!(E,\nabla_{\!E})=\pi_!(E',\nabla_{\!E'})+\pi_!(E'',\nabla_{\!''})\in K^0_{\text{flat}}(B)$, which fits with relation \eqref{suitexacteplate}.

Consider $E'$ as a subbundle of $E$.
The vertical exterior differential operator $d^{\nabla_{\!E}}$
respects the subbundle (over $B$) $\Omega(Z,E')$ of the vertical de Rham complex $(\Omega(Z,E),d^{\nabla_{\!E}})$ (remember that only the retriction of $\nabla_{\!E}$ to $Z$ is used to construct $d^{\nabla_{\!E}}$);
this filtration $0\subset\
\Omega(M,E')\subset\Omega(M,E)$ gives rise to some spectral sequence, and to some filtration $0\subset FH^\bullet(Z,E)\subset H^\bullet(Z,E)$
of the fibral cohomology of $E$.

The $(E_0,d_0)$-term of this spectral sequence is the direct sum of the fibral de Rham complexes of $E'$ and of $E''$; consequently, the $E_1$-term is the direct sum $\pi_!E'\oplus\pi_!E''$ of the fibral cohomology of $E'$ and of $E''$.

Let $s\colon E\to E'$ be a smooth vector bundle morphism such that $s\circ i$ is the identity of $E'$, then
$E''$ will
be identified with the subbundle ${\text{Ker}}s$ of $E$ so that $E$
will be identified with $E'\oplus E''$. Thus $E$ inherits two flat
connections $\nabla_{\!E}$ and $\nabla_{\!E'}\oplus\nabla_{\!E''}$, whose difference
is (as was used in lemma \ref{L:exactseqtransgr} in a nonflat context) a one form $\omega$ with values in ${\text{Hom}}(E'',E')$.
On any closed $E''$-valued form, $d^{\nabla_{\!E}}$ applies as $\omega\wedge$ so that
the operator $d_1$ of the spectral sequence is given by
\begin{equation}\label{E:above}
d_1=[\omega\wedge]\colon H^\bullet(Z,E'')\longrightarrow H^{\bullet+1}(Z,E')
\end{equation} 

It is elementary to check that this is exactly the linking maps of the exact sequence \eqref{diagrexplat}, and that the
spectral sequence converges at $E_2$
which is the filtrated fibral cohomology of $E$.

The exact diagram \eqref{diagrexplat} provides
(as in \eqref{complexestousrondsavant} and \eqref{E:sixdiagram}) a class of link $[\ell_{\rm{flat}}]$ between
$\pi_!^+E-\pi_!^-E$ and $(\pi_!^+E'\oplus\pi_!^+E'')-(\pi_!^-E'\oplus\pi_!^-E'')$,
the previous consideration establish that this class of link $[\ell_{\rm{flat}}]$ is exactly the same as the one obtained
by the construction following \eqref{E:linkcohom} from the
complex associated to $[\omega\wedge]$:
\begin{equation}\label{poiyu}\ldots\overset{[\omega\wedge]}\longrightarrow H^i(M,E')\oplus H^i(M,E'')\overset{[\omega\wedge]}\longrightarrow H^{i+1}(M,E')
\oplus H^{i+1}(M,E'')\overset{[\omega\wedge]}\longrightarrow\ldots\end{equation}
(which is the direct sum (with suitably shifted degrees) of maps \eqref{E:above} above) modulo the canonical isotopy class of isomorphism between (graded) cohomology and (graded) filtrated cohomology.
\begin{equation}\label{ralbol}
H^\bullet(Z,E)\cong FH^\bullet(Z,E)\oplus\big(H^\bullet(Z,E)/FH^\bullet(Z,E)\big)
\end{equation}
This is because the kernel of the map \eqref{E:above} is $H^i(M,E)/FH^i(M,E)$ and its cokernel is $FH^{i+1}(M,E)$, and because one can decompose \eqref{diagrexplat} in a direct sum of two
associated exact sequences (of length 4).
\subsubsection{Compatibility of direct images for $K^0_{\rm{top}}$ and $K^0_{\rm{flat}}$:}
By fibral Hodge theory, the $H^\pm(Z_y,E)$ are isomorphic to ${\rm{Ker}}
\big(d^{\nabla_{\!E}}+(d^{\nabla_{\!E}})^*\big)^\pm$ on $Z_y$ which are of constant dimension (whatever the riemannian metric on $M$ and the hermitian metric on $E$ may be), so that $(\{0\},\{0\},0)$ are suitable data in this situation. The compatibility of direct images in topological and
flat $K$-theories with the forgetful map $K_0^{\text{flat}}\to K_0^{\text{top}}$
trivially follows.

This can however be precised as following. The exact sequence
\eqref{suitexplate} provides the smooth isomorphism class $s\oplus p\colon E\longrightarrow E'\oplus E''$ between
$E$ and $E'\oplus E''$. Thus of course
\[ [\pi_!^+E]-[\pi_!^-E]=[\pi_!^+E'\oplus\pi_!^+E'']-
[\pi_!^-E'\oplus\pi_!^-E'']\quad{\text{ in }}\quad K^0_{\rm{top}}(B)\]
But as they are both obtained from the direct image constructions of the preceding section, one obtains some canonical equivalence class of link
$[\ell_{\rm{top}}]$ between $\pi_!^+E-\pi_!^-E$ and $(\pi_!^+E'\oplus\pi_!^+E'')-(\pi_!^-E'\oplus\pi_!^-E'')$
(see the penultimate alinea of \S\ref{clink}).
\begin{lemma}\label{difficile} $\ \ [\ell_{\rm{top}}]=[\ell_{\rm{flat}}]$  
\end{lemma}
\begin{proof}

For $\theta\in[0,1]$, put $\nabla_{\!\theta}=(\nabla_{\!E'}\oplus\nabla_{\!E''})
+\theta\omega$, then $\nabla_{\!E}=\nabla_{\!1}$, and $\nabla_{\!\theta}$ is flat
for any $\theta\in[0,1]$. Moreover, the flat bundles $(E,\nabla_{\!\theta})$
and $(E,\nabla_{\!E})$ are isomorphic for any $\theta>0$ throw the
automorphism ${\text{Id}}_{E'}\oplus\theta{\text{Id}}_{E''}$ of $E$.
For any $\theta$, $d^{\nabla_{\!\theta}}$ (as $d^{\nabla_{\!E}}$) also respects the subbundle $\Omega(Z,E')$ of $\Omega(Z,E)$, and the associated spectral sequence is isomorphic to the preceding one if $\theta$ is positive, so that the considerations of the preceding paragraph
apply verbatim for $\theta\in(0,1]$.

Put any riemanian metric on $M$, and endow $E\cong E'\oplus E''$ with a
direct sum hermitian metric. Hodge theory provides for any $\theta\in(0,1]$ an isomorphism between the (graded) kernel ${\mathcal H}_\theta^\bullet$ of the
fibral Dirac operator ${\mathcal D}_\theta=d^{\nabla_{\!\theta}}+(d^{\nabla_{\!\theta}})^*$
and the
cohomology of the de Rham complex associated
with $d^{\nabla_{\!\theta}}$.
In particular, the dimension of ${\mathcal H}^{i}_\theta$
is constant for any $i$ when $\theta$ goes over $(0,1]$. The isomorphism class provided by parallel
transport along $(0,1]$ of ${\mathcal H}_\theta$ is isotopic to the isomorphism provided
by twisting the de Rham cohomology by ${\text{Id}}_{E'}\oplus\theta{\text{Id}}_{E''}$.

Let $d^{\nabla_{\!E'}}$ and $d^{\nabla_{\!E''}}$ be the fibral exterior differential operators on
$\Omega(Z,E')$ and $\Omega(Z,E'')$ respectively which are obtained using $\nabla_{\!E'}$ and $\nabla_{\!E''}$, and 
define ${\mathcal D}'=d^{\nabla_{\!E'}}+
(d^{\nabla_{\!E'}})^*$ and ${\mathcal D}''=d^{\nabla_{\!E''}}+(d^{\nabla_{\!E''}})^*$. Then
${\mathcal D}_\theta={\mathcal D}'+{\mathcal D}''+\theta(\omega+\omega^*)$ so that one has a
continuous family of elliptic operators on $B\times[0,1]$.
Let $\lambda_{\text{min}}$ be the minimum positive eigenvalue of ${\mathcal D}'+{\mathcal D}''$ along all $B$.
(The compactness of $B$ ensures the existence (and the positivity!) of $\lambda_{\rm{min}}$).
There exists $\varepsilon>0$ such that $\theta\omega$ is bounded by $\frac15\lambda_{\rm{min}}$
in $L^2$ norm for all $\theta\leq\varepsilon$. Then for any $y\in B$ and any $\theta
\leq\varepsilon$, ${\mathcal D}_\theta$ has no eigenvalue
equal to $\pm\frac\lambda2$. Thus the (graded) direct sum ${\mathcal F}^\bullet_\theta$ of eigenspaces of
${\mathcal D}_\theta$ corresponding to eigenvalues belonging
to $]-\frac\lambda2,\frac\lambda2[$ is a finite rank vector bundle on $B\times
[0,\varepsilon]$ whose restriction ${\mathcal F}_0$ to $B\times\{0\}$ equals ${\rm{Ker}}{\mathcal D}'\oplus{\rm{Ker}}{\mathcal D}''$.

For any $\theta\in(0,\varepsilon]$, $d^{\nabla_{\!\theta}}$ respects ${\mathcal F}_\theta$ and
${\mathcal H}_\theta$ is identified with the cohomology
of the restriction of $d^{\nabla_{\!\theta}}$ to ${\mathcal F}_\theta$. This provides a link between ${\mathcal H}^+_\theta
-{\mathcal H}^-_\theta$ and ${\mathcal F}^+_\theta-{\mathcal F}^-_\theta$ by the construction of \eqref{E:linkcohom} \eqref{linkcohomprime}.
As was remarked in \eqref{oscour}, this link is nothing but the link obtained as in \eqref{precedente} from the exact sequence:
\[0\longrightarrow{\mathcal H}^+_\theta\longrightarrow{\mathcal F}^+_\theta
\overset{-{\mathcal D}_\theta}{-\!\!\!-\!\!\!-\!\!\!\longrightarrow}{\mathcal F}^-_\theta\longrightarrow{\mathcal H}^-_\theta
\longrightarrow0\]
And as was proved in the last alinea of \S\ref{clink}, this link is exactly the canonical link between ${\mathcal H}^+_\theta-{\mathcal H}^-_\theta$ and ${\mathcal F}^+_\theta-{\mathcal F}^-_\theta$ as different realisations of $\pi^{\rm{Eu}}_*[E]$. (The simplified version of \eqref{linkvaleurspropres} above
is due to the fact that ${\mathcal D}_\theta$
vanishes here on ${\mathcal H}_\theta^\pm$ for any $\theta$).

The composition of the Hodge isomorphism $\pi_!E\cong{\mathcal H}_1$, the parallel transport ${\mathcal H}_1\cong{\mathcal H}
_\theta$, the above link between ${\mathcal H}_\theta^+-{\mathcal H}^-_\theta$ and ${\mathcal F}_\theta^+-{\mathcal F}_\theta^-$, the parallel transport again ${\mathcal F}_\theta\cong{\mathcal F}_0$ and the Hodge isomorphism again ${\mathcal F}_0\cong\pi_!E'\oplus\pi_!E''$ provides a link between
$\pi_!^+E-\pi_!^-E$ and $(\pi_!^+E'\oplus\pi_!^+E'')-(\pi^-_!E'\oplus\pi^-_!E'')$. Its equivalence class is obviously independent of $\theta\in(0,\varepsilon]$
because of the compatibility of links between ${\mathcal F}_\theta$
and ${\mathcal H}_\theta$ with parallel transports for ${\mathcal F}$ and ${\mathcal H}$. This class of link is exactly $[\ell_{\rm{top}}]$, because of the fact that canonical links between different representatives of the topological direct image are mutually compatible.

It is straightforward to check that
$\frac1\theta P^{{\mathcal F}_\theta}d^{\nabla_{\!\theta}}P^{{\mathcal
F}_\theta}$ converges to $P^{{\mathcal F}_0}(\omega\wedge)P^{{\mathcal F}_0}$
as $\theta$ converges to $0$,
and this is the image of $[\omega\wedge]$ throw the Hodge isomorphism
${\mathcal F}_0\cong\pi_!E'\oplus\pi_!E''$.

Firstly, this proves that ${\mathcal H}_\theta$ converges to the kernel ${\mathcal H}_0$ of $P^{{\mathcal F}_0}(\omega\wedge)P^{{\mathcal F}_0}$ as
$\theta$ converges to $0$, because the dilation factor $\frac1\theta$ does not modify the kernels. This limit subspace ${\mathcal H}_0$ is identified by Hodge isomorphism ${\mathcal F}_0\cong H(Z,E')\oplus H(Z,E'')$ with the filtrated fibral cohomology of $E$ as seen around equation \eqref{E:above}, and this proves that
the parallel transport along $[0,1]$ for ${\mathcal H}$ provides the natural isotopy class of isomorphism between the fibral cohomology of $E$ and its filtrated counterpart \eqref{ralbol}.

Secondly, this proves that the link between ${\mathcal H}^+_\theta-{\mathcal H}^-_\theta$ and ${\mathcal F}^+_\theta
-{\mathcal F}^-_\theta$ converges to the link between
${\mathcal H}^+_0-{\mathcal H}^-_0$ and ${\mathcal F}^+_0-{\mathcal F}^-_0$ 
provided from \eqref{poiyu} by the construction \eqref{E:linkcohom} \eqref{linkcohomprime} (the finite rank of ${\mathcal F}_\theta$ makes these convergences elementary). The composition of these two last links is exactly $[\ell_{\rm{flat}}]$, and the independence on $\theta$ of the description of
$[\ell_{\rm{top}}]$ above implies the lemma.
\end{proof}
\subsubsection{Direct image of adjoint transpose flat bundles:}
Choose any ${\mathcal C}^\infty$ supplementary subbundle $T^H\!M$ of
$TZ$ in $TM$. Of course $T^H\!M\cong\pi^*TB$.
For any $y\in B$ and any vector ${\tt u}\in T_yB$, consider its horizontal lift ${\tt u}^H$:
${\tt u}^H$ is a global section of the restriction of $T^H\!M$ to $Z_y=\pi^{-1}(y)$ such that at
any point of $Z_y$ one has $\pi_*{\tt u}^H={\tt u}$.

Consider some vector bundle $\xi$ on $M$ with a (nonnecessarily flat)
connection $\nabla_{\!\xi}$ and hermitian metric $h^\xi$.
Remember the definition of ${\mathcal E}$ from \eqref{E:infinitedimvectbundle}. The flow associated to vector fields of the form ${\tt u}^H$ send fibres of $\pi$ to fibres of $\pi$ diffeomorphically, so that there is some fiberwise Lie differentiation operator ${\mathcal L}^{\nabla_{\!\xi}}_{{\tt u}^H}$ which acts on
$\xi$-valued vertical differential forms ${\mathcal E}$ (it is of course defined using the connection $\nabla_{\!\xi}$).
Put then for any local section $\sigma$ of ${\mathcal E}$
\begin{equation}\label{E:megaconnexion}
\overline\nabla_{\tt u}\sigma={\mathcal L}^{\nabla_{\!\xi}}_{{\tt u}^H}\sigma
\end{equation}
(see \cite{BismutLott} definition 3.2). $\overline\nabla$ is a connection on ${\mathcal E}$
as can be proved following \cite{BismutLott} (3.8) to (3.10).

Define the adjoint connection associated to $\overline\nabla$ as in paragraph \ref{adjcon} by the following formula, valid
for any element ${\tt u}$ of the tangent bundle of $B$ and any local sections
$\sigma$ and $\theta$ of ${\mathcal E}$:
\begin{equation}\label{E:connexions}
\langle\overline\nabla_{\tt u}^S\sigma,\theta\rangle_{L^2}={\tt u}.\langle \sigma,\theta\rangle_{L^2}-\langle
\sigma,\overline\nabla_{\tt u} \theta\rangle_{L^2}
\end{equation}

Consider the adjoint transpose $\nabla_{\!\xi}^*$ of $\nabla_{\!\xi}$, and
the associated connection $\overline\nabla\check{\ }$ on ${\mathcal E}$
(in the same way as $\overline\nabla$ is associated to $\nabla_{\!\xi}$ throw \eqref{E:megaconnexion}), call $\overline\nabla\check{\ }^S$ the
adjoint connection associated to $\overline\nabla\check{\ }$ as in \eqref{E:connexions}.
The reader is warned that the connection denoted by $\overline\nabla^S$
corresponds to the connection denoted by $\widetilde\nabla^*$ in \cite{BismutLott} Proposition 3.7 and that $\overline\nabla\check{\ }$ and $\overline\nabla\check{\ }^S$
here have no counterpart in \cite{BismutLott}.
\begin{lemma}\label{commu*}For any vector ${\tt u}$ tangent to $B$, and any local section $\sigma$ of ${\mathcal E}$
\[\overline\nabla_{\tt u}\!\!\check{\ }=*_Z^{-1}(\overline\nabla_{\tt u}^S(*_Z\sigma))\qquad{\text{and}}\qquad
\overline\nabla_{\tt u}\!\!\check{\ }^{S}
=*_Z^{-1}(\overline\nabla_{\tt u}(*_Z\sigma))\]
\end{lemma}
\begin{proof} Let $P^{TZ}$ be the projection of $TM$ onto $TZ$ with
kernel $T^H\!M$.
All riemannian metrics on $M$ which coincide with $g^Z$
on $TZ$ and make $TZ$ and $T^H\!M$ orthogonal give rise to
Levi-Civita connections $\nabla_{\!LC}$ on $TM$ which all
project to the same connection $\nabla_{\!TZ}=P^{TZ}\nabla_{\!LC}$
on $TZ$. Denote allways by $\nabla_{\!TZ}$ the associated
connection on $\wedge^\bullet T^*Z$, it is compatible with the
Clifford action \eqref{Clifford}, so that its associated covariant derivative commutes with $*_Z$.
Let $\nabla_{\!TZ\otimes\xi}$ be the connection
on $\wedge^\bullet T^*Z\otimes\xi$ associated to $\nabla_{\!TZ}$ and $\nabla_{\!\xi}$, its transpose adjoint $\nabla_{\!TZ\otimes\xi}^*$ with respect to $(\ /\ )_Z$
is nothing but the connection on $\wedge^\bullet T^*Z\otimes\xi$
associated to $\nabla_{\!TZ}$ and $\nabla_{\!\xi}^*$.
Then the covariant derivatives associated to both $\nabla_{\!TZ\otimes\xi}$ and $\nabla_{\!TZ\otimes\xi}^*$
commute with $*_Z$.

Let ${\tt u}$ be some vector tangent to $B$, and ${\tt u}^H$ its horizontal lift.
For any vertical vector ${\tt y}$ tangent to te fibre, the vertical projection
$P^{TZ}{\nabla_{\!LC}}_{\tt y}^{\ }{\tt u}^H$ of
the covariant derivative throw the connection $\nabla_{\!LC}$ of ${\tt u}^H$ along ${\tt y}$ is independent of the
global riemannian metric defining $\nabla_{\!LC}$. Moreover,
if ${\tt v}$ is another vertical tangent vector at the same point as ${\tt y}$, then the scalar product
$g^Z(P^{TZ}{\nabla_{\!LC}}_{\tt y}^{\ }{\tt u}^H,{\tt v})$ is symmetric in ${\tt y}$ and ${\tt v}$.

As proved in \cite{BismutLott} (3.27) and (3.32), if $({\tt e}_1,{\tt e}_2,\ldots,{\tt e}_{{\rm{dim}}Z})$ is
an orthonormal base of $TZ$, then for any local section $\sigma$ of ${\mathcal E}$, the
connections $\overline\nabla$ and $\overline\nabla^S$ express locally on $M$ as:
\begin{equation}
\begin{aligned}
\overline\nabla_{\tt u}^{\ }\sigma&={\nabla_{\!TZ\otimes\xi}}_{{\tt u}^H}^{\ }\sigma+\sum_{i{\text{ and }}k}
g^Z(P^{TZ}{\nabla_{\!LC}}_{{\tt e}_i}^{\ }{\tt u}^H,{\tt e}_k){\tt e}_i^\flat\wedge(\iota_{{\tt e}_k}\sigma)\\
\overline\nabla^S_{\tt u}\sigma&={\nabla_{\!TZ\otimes\xi}^*}_{{\tt u}^H}\sigma-\sum_{i{\text{ and }}k}
g^Z(P^{TZ}{\nabla_{\!LC}}_{{\tt e}_i}^{\ }{\tt u}^H,{\tt e}_k){\tt e}_i^\flat\wedge(\iota_{{\tt e}_k}\sigma)
\end{aligned}
\end{equation}
The lemma follows from the obvious corresponding formulae for $\overline\nabla\check{\ }$ and $\overline\nabla\check{\ }^{S}$, the fact that $\nabla_{\!TZ\otimes\xi}$
and $\nabla_{\!TZ\otimes\xi}^*$ commute with $*_Z$, the symmetry in ${\tt e}_i$ and ${\tt e}_k$ of
$g^Z(P^{TZ}{\nabla_{\!LC}}_{{\tt e}_i}^{\ }{\tt u}^H,{\tt e}_k)$ and the fact that
for any $i$ and $k$
\[({\tt e}_i^\flat\wedge)\iota_{{\tt e}_k}*_Z=-*_Z({\tt e}_k^\flat\wedge)\iota_{{\tt e}_i}\]

Another (yet nonlocal) proof of this result can be obtained from the fact that the fiberwise Lie differentiation verifies a formula like \eqref{flipflop}
\[{\mathcal L}_{\tt u}(\sigma\underset{h^\xi}\wedge\theta)=
({\mathcal L}^{\nabla_{\!\xi}}_{\tt u}\sigma)\underset{h^\xi}\wedge\theta+
\sigma\underset{h^\xi}\wedge({\mathcal L}^{\nabla^*_{\!\xi}}_{\tt u}\theta)\]
where ${\mathcal L}$ is the usual fiberwise Lie derivative of ordinary vertical differential forms. Applying this to vector fields of the form ${\tt u}^H$ and integrating along the fibres of $\pi$ yields
\[{\tt u}.\left(\int_Z\sigma\underset{h^\xi}\wedge\theta\right)=
\int_Z(\overline\nabla_{\tt u}\sigma)\underset{h^\xi}\wedge\theta+
\int_Z\sigma\underset{h^\xi}\wedge(\overline\nabla_{\tt u}\!\!\check{\ }\theta)\]
From this last formula and \eqref{L2glob}, \eqref{starwedge} and \eqref{E:connexions}, the lemma may be proved by reasoning as in
the first part of the proof of \cite{BismutLott} Theorem 2.32.
\end{proof}
Back to flat bundles, suppose that $E$ is a vector bundle on $M$ with a flat connection $\nabla_{\!E}$ (and hermitian metric $h^E$), and construct the associated objects ${\mathcal E}$,
$\overline\nabla$ and $\overline\nabla^S$ as above.
Let $P\colon{\mathcal E}\longrightarrow{\rm{Ker}}\big(d^{\nabla_{\!E}}+
(d^{\nabla_{\!E}})^*\big)$ be the orthogonal projection, then it is proved in \cite{BismutLott} Proposition 3.14 that $\nabla_{\pi_!E}\cong P \overline\nabla P$ and $\nabla_{\pi_!E}^*\cong P\overline\nabla^SP$
throw the fibral Hodge isomorphism $\pi_!E\cong{\rm{Ker}}\big(d^{\nabla_{\!E}}+
(d^{\nabla_{\!E}})^*\big)$.

Consider on $E$ the adjoint transpose connection $\nabla_{\!E}^*$
(which is flat). The direct image of the flat bundle $(E,\nabla_{\!E}^*)$
will be denoted by $\pi_!\!{\ }E$ and the flat connection on it by $\nabla_{\!\pi_!\!\check{\ }E}$ so that $\pi_!(E,\nabla_{\!E}^*)
=(\pi_!\!\check{\ }E,\nabla_{\!\pi_!\!\check{\ }E})$.

As precedingly, call $\overline\nabla\check{\ }$ and $\overline\nabla\check{\ }^S$
the connections on ${\mathcal E}$ constructed from $\nabla_{\!E}^*$ as in \eqref{E:megaconnexion} and \eqref{E:connexions}, let $P\check{\ }\colon{\mathcal E}\longrightarrow{\rm{Ker}}\big(d^{\nabla_{\!E}^*}+
(d^{\nabla^*_{\!E}})^*\big)$ be the orthogonal projection, then from \cite{BismutLott} Proposition 3.14 again, $\nabla_{\!\pi_!\!\check{\ }E}\cong P\check{\ } \overline\nabla\check{\ } P\check{\ }$ and $\nabla_{\!\pi_!\!\check{\ }E}^*\cong P\check{\ }\overline\nabla\check{\ }^SP\check{\ }$
throw the fibral Hodge isomorphism $\pi_!\!\check{\ }E\cong{\rm{Ker}}\big(d^{\nabla_{\!E}^*}+
(d^{\nabla_{\!E}^*})^*\big)$.

It follows from \eqref{agentdechange} that $P\check{\ }=*_ZP*_Z^{-1}$ so that
$*_Z$ directly provides a smooth isomorphism $\pi_!E\cong\pi_!\!\check{\ }E$ (which is of course exactly $[\ell_*]$ of \S\ref{evenHodge} if ${\rm{dim}}Z$ is even, and
compatible with the constructions of $[\ell^{\{0\}}_{{\mathcal H}_0}]$ and $[\ell^{\{0\}}_{{\mathcal H}_1}]$ in \S\ref{oddHodge} if ${\rm{dim}}Z$ is odd).
It then follows
from the preceding lemma \ref{commu*} that throw this isomorphism
$\nabla_{\!\pi_!E}^*\cong\nabla_{\!\pi_!\!\check{\ }E}$ and $\nabla_{\!\pi_!E}\cong
\nabla_{\!\pi_!\!\check{\ }E}^*$. Now $*_Z$ respects the ${}^+$ and ${}^-$ parts of
${\mathcal E}$ if ${\rm{dim}}Z$ is even, and exchanges them if ${\rm{dim}}Z$ is odd so that
\begin{lemma}\label{dirimconjflat} The following equalities hold in $K^0_{\rm{flat}}(B)$:
\begin{align*}
\pi_!(E,\nabla_{\!E}^*)&=(\pi_!E,\nabla_{\!\pi_!E}^*)\qquad{\text{ if ${\rm{dim}}Z$ is even}}\\
\pi_!(E,\nabla_{\!E}^*)&=-(\pi_!E,\nabla_{\!\pi_!E}^*)\qquad{\text{ if ${\rm{dim}}Z$ is odd}}
\end{align*}
\end{lemma}
This means that $\pi_!$ commutes or anticommutes with the conjugation in $K^0_{\rm{flat}}$ of \eqref{collier} if ${\rm{dim}}Z$ is even or odd
respectively.
In particular, $\pi_!$ respects or exchanges the ``real'' and ``imaginary''
parts of $K^0_{\rm{flat}}$ (with respect to \eqref{collier}) if
${\rm{dim}}Z$ is even or odd respectively.
\begin{remark} If ${\rm{dim}}Z$ is odd, let $[\ell^{\{0\}}_{\pi_!E}]$ be the link between
$\pi_!^+E-\pi_!^-E$ and $\{0\}-\{0\}$ of lemma \ref{linkzero}, then the map \[(E,\nabla_{\!E})\in K^0_{\rm{flat}}(M)\longmapsto
\big(\pi_!^-E,\nabla_{\!\pi_!^-E},\pi_!^+E,\nabla_{\!\pi_!^+E},[-\ell^{\{0\}}_{\pi_!E}]\big)
\in K^0_{\rm{rel}}(B)\] 
defines a group morphism $\pi_\leftarrow\colon K^0_{\rm{flat}}(M)\to K^0_{\rm{rel}}(B)$
such that $\pi_!=\partial\circ\pi_\leftarrow$.
This is because of lemma \ref{difficile} and of the obvious compatibility of the link of lemma \ref{linkzero} with direct sums.
It follows from lemma \ref{dirimconjflat} above that
$\pi_\leftarrow$ anticommutes with the conjugations \eqref{collier} on $K^0_{\rm{flat}}(M)$
and \eqref{congrel} on $K^0_{\rm{rel}}(B)$.
\end{remark}
\subsection{Relative $K$-theory:}\label{dirimrelsection}

For some element $(E,\nabla_{\!E},F,\nabla_{\!F},f)$
of $K_0^{\text{rel}}(M)$,
call $\nabla_{\!\pi_!E}^+$, $\nabla_{\!\pi_!E}^-$, $\nabla_{\!\pi_!F}^+$ and $\nabla_{\!\pi_!F}^-$ the respective flat connections on the bundles $\pi_!^+E$, $\pi_!^-E$, $\pi_!^+F$ and $\pi_!^-F$ on $B$, (obtained from the direct image of the
flat bundles $(E,\nabla_{\!E})$ and $(F,\nabla_{\!F})$), then the construction
of the penultimate alinea of \S\ref{clink} produces a link $[\ell_{\pi_*f}]$
between $\pi_!^+E-\pi_!^-E$ and $\pi_!^+F-\pi_!^-F$.
\begin{definition}\label{dirimrel}
\begin{align*}
\pi_*(E,\nabla_{\!E},F,\nabla_{\!F},f)=\big(
&
\pi^+_!E\oplus\pi_!^-F,\nabla_{\!\pi^+_!E}\oplus\nabla_{\!\pi^-_!F},
\\&\qquad 
\pi_!^-E\oplus\pi^+_!F,\nabla_{\!\pi^-_!E}\oplus\nabla_{\!\pi_!^+F},
[\ell_{\pi_*f}]\big)
\end{align*}
\end{definition}
\begin{theorem*}
This defines a morphism $K_0^{\text{rel}}(M)
\to K_0^{\text{rel}}(B)$ which enters in the following commutative diagram (with lines modeled on \eqref{suitexmoi}):
\[\begin{CD}
K^1_{\text{top}}(M)@>>>K^0_{\text{rel}}(M)
@>\partial>>K^0_{\text{flat}}(M)@>>>K^0_{\text{top}}(M)\\
@V\pi^{\rm{Eu}}_*VV@V\pi_*VV@V\pi_!VV@VV\pi^{\rm{Eu}}_*V\\
K^1_{\text{top}}(B)@>>>K^0_{\text{rel}}(B)
@>\partial>>K^0_{\text{flat}}(B)@>>>K^0_{\text{top}}(B)
\end{CD}\]
If ${\rm{dim}}Z$ is odd, then $\pi_*=\pi_\leftarrow\circ\partial$.
\end{theorem*}
\begin{proof}
The compatibility of the construction of $\pi_*$
with direct sums is tautological, it is also clear that $\pi_*(E,\nabla_{\!E},F,\nabla_{\!F},f)$ depends
on $f$ only through its isotopy class. So relations $(i)$
and $(ii)$ in the definition of $K^0_{\text{rel}}(M)$
are sent to relations in $K^0_{\text{rel}}(B)$, and lemma \ref{difficile} is exactly what is needed here to prove that relation $(iii)$ in the definition of $K^0_{\rm{rel}}(M)$ is sent to a relation of the same type in $K^0_{\rm{rel}}(B)$.

The commutativity of the diagram should be clear if one remembers that the direct image $\pi_*^{\rm{Eu}}\colon K^1_{\text{top}}(M)\to K^1_{\text{top}}(B)$ will
simply be obtained from the direct image $\big(\pi\times{\rm{Id}}_{S^1}\big)^{\rm{Eu}}_*\colon K^0_{\text{top}}
(M\times S^1)\to K^0_{\text{top}}(B\times S^1)$, this construction is obviously compatible with the one used in the construction of $[\ell_{\pi_*f}]$. 

It follows from lemma \ref{linkzero} that if $(E,\nabla_{\!E},F,\nabla_{\!F},f)\in K^0_{\rm{rel}}(M)$ and if ${\rm{dim}}Z$ is odd and, then
$[\ell^{\{0\}}_{\pi_!F}]=[\ell_{\pi_*f}]\circ[\ell^{\{0\}}_{\pi_!E}]$ (where
$[\ell^{\{0\}}_{\pi_!F}]$ or $[\ell^{\{0\}}_{\pi_!E}]$ is the link between
$\pi_!^+E-\pi_!^-E$ or $\pi_!^+F-\pi_!^-F$ respectively and $\{0\}-\{0\}$ of
lemma \ref{linkzero}) so that $[\ell_{\pi_*f}]=[\ell^{\{0\}}_{\pi_!F}]\circ[\ell^{\{0\}}_{\pi_!E}]^{-1}$ and
\begin{align*}\pi_*(E,\nabla_{\!E}, F,\nabla_{\!F},f)&=
(\pi_!F^+,\nabla_{\pi_!F^+},\pi_!F^-,
\nabla_{\pi_!F^-},[\ell^{\{0\}}_{\pi_!F}])\\&\qquad-
(\pi_!E^+,\nabla_{\pi_!E^+},\pi_!E^-,
\nabla_{\pi_!E^-},[\ell_{\pi_!E}^{\{0\}}])
\end{align*}
which proves the last statement of the theorem.\end{proof}
If ${\rm{dim}}Z$ is even, it follows from lemmas \ref{starlink} and \ref{dirimconjflat} that
\begin{equation}\label{dirimsymmetry}\pi_*(E,\nabla_{\!E},E,\nabla_{\!E}^*,{\rm{Id}}_E)=(\pi_!E,\nabla_{\pi_!E},\pi_!E,\nabla_{\!\pi_!E}^*,{\rm{Id}}_{\pi_!E})\end{equation}

\end{document}